\documentclass[11pt]{article}

\usepackage[cp850]{inputenc}
\usepackage{latexsym,graphicx}

\font\tenmath=msbm10 scaled 1200
\font\sevenmath=msbm7 scaled 1200
\font\fivemath=msbm5 scaled 1200

\newfam\mathfam \textfont\mathfam=\tenmath
\scriptfont\mathfam=\sevenmath \scriptscriptfont\mathfam=\fivemath
\def\math{\fam\mathfam}
\def\R{{\math R}}
\def\N{{\math N}}
\def\E{{\math E}}

\def\P{{\math P}}

\def\Q{{\math Q}}

\def \^#1{\if#1i{\accent"5E\i}\else{\accent"5E#1}\fi}

\def \ni{\noindent}

\newtheorem{Thm}{Theorem}

\newtheorem{Lem}{Lemma}
\newtheorem{Pro}{Proposition}
\newtheorem{Cor}{Corollary}
\newtheorem{Dfn}{Definition}

\author{
{\sc Siegfried Graf} \thanks{Universit\"at Passau, Fakult\"at f\"ur Mathematik
und Informatik, D-94030
Passau, Germany. E-mail: {\tt  graf@fmi.uni-passau.de}}  \quad
{\sc Harald Luschgy}\thanks{Universit\"at Trier, FB IV-Mathematik, D-54286
Trier, Germany.
E-mail: {\tt luschgy@uni-trier.de}} \quad {\sc and}
\quad {\sc Gilles Pag\`es} \thanks{Laboratoire de Probabilit\'es et Mod\`eles
al\'eatoires, UMR~7599, Universit\'e Paris 6, case 188, 4,
pl. Jussieu, F-75252 Paris Cedex 5, France. E-mail:{\tt  gpa@ccr.jussieu.fr}}
}

\title{\bf Optimal quantizers for Radon random vectors in a Banach space}

\begin{document}


\maketitle
\begin{abstract}
For $n \in \N, r \in (0, \infty)$ and a Radon random vector $X$ with values in
a Banach space $E$ let $e_{n,r}(X,E) = \inf ( \E \min_{a \in \alpha} \| X-a
\|^r )^{1/r}$, where the infimum is taken over all subsets $\alpha$ of $E$
with card$(\alpha) \leq n$ ($n$-quantizers). We investigate the existence of
optimal $n$-quantizers for this $L^r$-quantization propblem, derive their
stationarity properties and establish for $L^p$-spaces $E$ the pathwise
regularity of stationary quantizers.
\end{abstract}


\bigskip
\ni {\em 2000 Mathematics Subject Classification: 41A46, 60B11, 94A29}

\section{Introduction}
\setcounter{equation}{0}
\setcounter{Assumption}{0}
\setcounter{Theorem}{0}
\setcounter{Proposition}{0}
\setcounter{Corollary}{0}
\setcounter{Lemma}{0}
\setcounter{Definition}{0}
\setcounter{Remark}{0}
We investigate optimal quantizers and the quantization error in the functional
$L^r$-quantization problem for stochastic processes viewed
as random variables in a Banach (function) space. So let $(E, \| \cdot \| )$
be
a real Banach space and consider a Radon random variable
$X : ( \Omega, {\cal A}, \P) \rightarrow E$ which means that $X$ is Borel
measurable and its
distribution $\P_X$ is a Radon probability measure on $E$. For $n \in \N$ and
$r \in (0, \infty)$, the $L^r$-quantization problem for $X$ of level $n$
consists in minimizing
\[
( \E \min\limits_{a \in \alpha} \|  X - a \|^r)^{1/r} = \|  \min_{a \in
\alpha} \| X - a \| \|_{L^r(\P)}
\]
over all subsets $\alpha \subset E$ with $\mbox{card} (\alpha) \leq n$. Such a
set $\alpha$ is
called $n$-codebook or $n$-quantizer. The minimal $n$th quantization error is
then defined by
\begin{equation}
e_{n,r} (X,E) := \inf \{ (\E \min_{a \in \alpha} \| X-a \|^r)^{1/r}  : \alpha
\subset E, \; \mbox{card} (\alpha) \leq n\} .
\end{equation}
Under the integrability condition
\begin{equation}
\E \| X \|^r< \infty
\end{equation}
the quantity $e_{n,r} (X, E)$ is finite.

For a given $n$-codebook $\alpha$ one defines an associated closest neighbour
projection
\[
\pi_{\alpha} := \sum\limits_{a \in \alpha} a 1_{C_{a}(\alpha)}
\]
and the induced $\alpha-quantized \; version$ (or $\alpha-quantization$) of
$X$ by
\begin{equation}
\hat{X}^{\alpha} := \pi_\alpha (X) ,
\end{equation}
where $\{C_a (\alpha) : a \in \alpha \}$ is a Voronoi partition induced by
$\alpha$, that is a
Borel partition of $E$ satisfying
\[
C_a (\alpha ) \subset \{ x \in E : \| x-a \| = \min\limits_{b \in \alpha} \|
x-b \|  \}
\]
for every $a \in \alpha$. Then one easily checks that, for any measurable
random variable
$X^{'} : \Omega \rightarrow \alpha \subset E$,
\[
\E \| X - X^{'} \|^r \geq \E \| X - \hat{X}^\alpha \|^r
= \E \min\limits_{a \in \alpha} \| X - a \|^r
\]
so that finally
\begin{eqnarray}
e_{n,r} (X,E) & = & \inf \{ ( \E \| X - \hat{X} \|^r )^{1/r} :
 \hat{X} = f(X) , f : E \rightarrow E \; \mbox{Borel measurable}, \mbox{card}
\; (f(E)) \leq n \} \nonumber \\
 & = & \inf \{ ( \E \| X - \hat{X} \|^r)^{1/r} : \hat{X} : \Omega
       \rightarrow E \; \mbox{measurable},\;\mbox{card} \; (\hat{X} ( \Omega))
\leq n \} .
\end{eqnarray}

Functional quantization of stochastic processes can thus be seen as a
discretization of the path-space $E$ of a
process and the approximation (coding) of a stochastic process by finitely
many deterministic functions from its path-space. Typical settings are $E =
L^p([0,1], dt)$ and $E = C([0,1])$. Functional quantization is the natural
extension to stochastic processes or Banach space valued random vectors of
the so-called optimal vector quantization of random vectors in $E = \R^d$
which has been extensively investigated since the late 1940's in Signal
processing and Information Theory (see \cite{GERSH}, \cite{GRAY}). For the
mathematical aspects of vector quantization in $\R^d$, one may consult
\cite{GRLU} and for algorithmic aspects see \cite{PAGEPRIN}.

Recently, the extension of optimal vector quantization to stochastic processes
has given raise to many theoretical
developments including the rate of convergence of the quantization errors
$e_{n,r} (X)$ to zero as $n \rightarrow \infty$ and the construction of
good or even rate optimal quantizers (see e.g. \cite{DEFEMASC}, \cite{DEREI}, \cite{DEREI1}
\cite{GRLUPA}, \cite{LUPA1}, \cite{LUPA2}, \cite{LUPA3}). For a first
promising application to the pricing of financial derivatives through
numerical integration on path-spaces see \cite{PAPR}. In this paper we aim to
develop general results on the existence of optimal quantizers and their
properties.

The paper is organized as follows. In Section 2, a theorem about the existence
of optimal $n$-quantizers for
$E$-valued Radon random vectors lying in $E$ or in some suitable superspace $G
\supset E$ is established under some very general assumptions. It
relates existence to intersection properties of closed balls. This problem is
connected with its bidual counterpart and enlightened by
counterexamples. Furthermore, bounds of the quantization errors $e_{n,r}(X,E)$
in terms of $e_{n,r}(X,G)$ for superspaces $G$ and in terms of
marginals of $X$ for vector valued processes are derived.
In Section 3 the stationarity property of optimal $n$-quantizers is
investigated. This turns out to be an essential key for
the functional quantization of 1-dimensional diffusion processes (see
\cite{LUPA3}). 
For smooth Banach spaces stationary quantizers are defined as the critical points of the 
distortion function. In the case of $L^p$-spaces $E$ which are natural path-spaces of processes
some pathwise regularity for these stationary quantizers is established.
The result applies e.g. to Gaussian processes,
$d$-dimensional diffusion processes and certain L\'evy processes. \\ \\
\section{Optimal quantizers and quantization errors}
%
\setcounter{equation}{0}
\setcounter{Assumption}{0}
\setcounter{Theorem}{0}
\setcounter{Proposition}{0}
\setcounter{Corollary}{0}
\setcounter{Lemma}{0}
\setcounter{Definition}{0}
\setcounter{Remark}{0}
Let $X$ be a Radon $(E, \| \cdot \| )$-valued random variable with
distribution $\P_X$.
The Radon property of $\P_X$ means inner regularity w.r.t. compact sets and on
Banach spaces it is the same as tightness which in turn is equivalent to the
existence of a separable Borel measurable set with $\P_X$-probability 1.
It is to be noticed that if $\P (X \in F) = 1$ for some Banach subspace $F$ of
$E,X$ is Radon when viewed as $F$-valued random variable. On the other hand,
if $E$ is a Banach
subspace of some Banach space $G$ then $X$ is also Radon as $G$-valued random
variable.

We will assume throughout this section that $X$ satisfies the integrability
condition (1.2) for some
$r \in (0, \infty)$. Then
\begin{equation}
\lim_{n \to \infty} e_{n,r} (X,E) = 0.
\end{equation}
As a matter of fact, the support of $\P_X$ being separable there exists a
countable subset $\{ a_n, n \geq 1 \}$ everywhere dense in
$\mbox{supp} (\P_X)$. It is clear that
\[
0 \leq e^r_{n,r} (X,E) \leq \E \min_{1 \leq i \leq n} \| X - a_i \|^r
\rightarrow 0 \; \mbox{as} \; n \rightarrow \infty
\]
by the Lebesgue dominated convergence Theorem. On the other hand, the
existence of optimal quantizers, i.e. the fact that
$e_{n,r} (X,E)$ actually stands as a minimum needs much more care. \\
\medskip
\subsection{Existence of optimal quantizers}
$A$ set $\alpha \subset E \; \mbox{with} \; 1 \leq \; \mbox{card} \; ( \alpha)
\leq n$
is called an {\em $L^r$-optimal $n$-quantizer} for $X$ if
\begin{equation}
(\E \min_{a \in \alpha} \| X-a \|^r)^{1/r} = e_{n,r} (X,E).
\end{equation}
Let ${\cal C}_{n,r} (X,E)$ denote the set of all $L^r$-optimal $n$-quantizers
for $X$ in $E$.

We first provide some interesting properties of $n$-optimal quantizers (they
can be seen
as necessary conditions for  $n$-optimality). Their proofs are literally the
same as those (established in  finite-dimension) of
Theorem 4.1 and Theorem 4.2 in \cite{GRLU} respectively. They are related with
the {\em Voronoi partitions} induced by a $n$-quantizer
$\alpha$: these are the Borel partitions $ \{ C_a (\alpha) : a \in \alpha\}$
of $E$  which satisfy
\begin{equation}
C_a(\alpha)\subset V_a(\alpha):=\left\{x \!\in E : \|x - a \|  = \min_{b \in
\alpha}
\| x - b \|
\right\}.
\end{equation}
Let us note that $V_a (\alpha)$ is closed and star-shaped relative for $a$ and
for every $a \in\alpha$,
\[
\left\{x \!\in E : \; \|x - a \|  < \min_{b \in \alpha \setminus \{ a \}}
\|x- a \|
\right\}\subset \,\stackrel{ \circ}{C}_a (\alpha) \subset
\overline{C_a(\alpha)} \subset V_a (\alpha).
\]
Furthermore,  as soon as $(E,\|\,.\,\|)$ is {\em strictly convex}
(\footnote{$i.e.$
$B_{E}(0,1)$ is a strictly convex set: $\forall x,y\!\in S_{E} (0,1),\,x\neq
y,\, \forall\lambda\!\!\in(0,1),\,
\|\lambda x+(1-\lambda)y\|<1$.}), any Voronoi partition   satisfies for every
$a \in\alpha$
\begin{equation}\label{adheranceC}
\overline{C_a(\alpha)}= V_a(\alpha)
\end{equation}
and
\[
\stackrel{ \circ}{C}_a(\alpha) = {\stackrel{\circ}{V}_a(\alpha)} = \left\{x
\!\in E : \; \|x - a \|  < \min_{b \in \alpha \setminus \{ a \}}
\| x - b \|
\right\}.
\]
\begin{Pro}
\label{Pro2.1&2} Assume that
card$(\mbox{\rm supp} (\P_{_X}) )\geq n$.

\smallskip
\noindent $(a)$ Let $\alpha \in {\cal C}_{n,r} (X,E)$. Then card$( \alpha ) =
n$ and for every $a \! \in \alpha$,
\[
\P_{_X}(C_a(\alpha)) > 0 \quad \mbox{ and } \quad \{ a \} \in {\cal C}_{1,r}
(\P_{_X}(\, \cdot\, | C_a (\alpha)), E) .
\]

\noindent $(b)$  Assume that $E$ is smooth {\rm (\footnote{i.e. the norm is
Gateaux-differentiable at every $x\neq 0 $.})} and strictly
convex.  If $\alpha\in {\cal C}_{n,r} (X,E)$ and
\[
(r > 1)
\qquad\mbox{ or }\qquad ( r = 1 \; \mbox{ and } \; \P(X\in \alpha) = 0) ,
\]
then
\begin{equation}\label{tessellation}
\P_{_X}(V_a(\alpha)\cap  V_b (\alpha)) = 0 \quad \mbox{ for every } \; a, b
\in \alpha, a \not= b.
\end{equation}
\end{Pro}

Note that, under the strict convexity assumption, (2.5) is then equivalent to
both
\[
\left(\forall\, a \in \alpha ,\quad \P_{_X}(\partial C_a(\alpha)) = 0\right)\;
\mbox{ and }\;
\left(\forall\, a \in \alpha ,\quad  \P_{_X}(\partial V_a(\alpha)) = 0\right).
\]

The first results of existence for optimal quantizers are due to
Cuesta-Albertos and Matr\`an \cite{CUEST} and P{\"a}rna \cite{PAERN}) for
uniformly convex
and reflexive Banach spaces, respectively. We provide an extension to Banach
spaces having the property that the closed balls form a
{\em compact system}.
A system ${\cal K}$ of subsets of $E$ is called compact if each subsystem
${\cal K}_0$ of
${\cal K}$ which has the finite intersection property $(i.e.$ the intersection
of each finite
subsystem of ${\cal K}_0$ is not empty) has a nonempty intersection. Let
$B(s,\rho) = B_E(x, \rho) := \{ y \in E: \| y-x \| \leq \rho \}$ be the closed
ball
of radius $\rho$ centered at $x$. \\
\begin{Dfn} A pair $(F,G)$ consisting of a Banach space $G$ and a Banach
subspace $F$ of
$G$ is called admissible if $\{B_G (x, \rho) : x \in F, \rho > 0 \}$ is a
compact system in $G$.
$G$ is called admissible if $(G,G)$ is admissible.
\end{Dfn}

The level $n$ $L^r$-distortion function is defined by
\begin{equation}
D^X_{n,r} : E^n \rightarrow \R_{+}, D^X_{n,r}(a) := \E \min_{1 \leq i \leq n}
\| X - a_i \|^r .
\end{equation}
\begin{Thm} Assume that $\P_X (F) = 1$ for some Banach subspace $F$ of $E$ and
that
$(F,E)$ is admissible. Then, for every $n \in \N$,
\[
{\cal C}_{n,r} (X,E) \not= \emptyset .
\]
\end{Thm}
\medskip
\noindent
{\bf Proof.} Fix $n \in \N$. Let $\tau_0$ denote the topology on $E$ generated
by the
system
$\{B(x,\rho)^c : x \!\in F,\,\rho > 0\}$ and let
$\tau$ be the product topology on $E^n$ (these topologies usually do not
satisfy the Hausdorff axiom). The
family $\{B(x,\rho) : x \!\in F,\,\rho > 0\}$ being a compact system in $E$,
one checks that $E$ is
$\tau_0$-quasi-compact(\footnote{$i.e.$ satisfies the Borel-Lebesgue axiom --
from any open covering one may extract a
finite open covering -- but possibly not the Hausdorff axiom.}). Consequently,
$E^n$ is $\tau$-quasi-compact. It is  obvious that any lower semi-continuous
(l.s.c.) function defined on $E^n$
then reaches  a minimum. Hence, the proof amounts to showing that the
distortion function
$D^X_{n,r}: E^n \rightarrow \R_+$
is $\tau$-lower semi-continuous.

For every $x \!\in F$ and  $a \in E^n$, set $\displaystyle d(x,a) := \min_{1
\leq i \leq n} \| x - a_i \|$. Then
\[
\left\{ a \in E^n :d (x , \cdot)^r \leq c \right\} = \bigcup\limits^n_{i=1} \{
a \in E^n : a_i \!\in B(x,c^{1/r}) \}
\]
is $\tau$-closed for every $c \geq 0$. Hence, $a\mapsto d(x, a)^r$ is
$\tau$-lower semi-continuous. In turn
any convex combination of such functions are $\tau$-l.s.c. as well. This
implies that $D^X_{n,r}$ (and $(D^X_{n,r})^{1/r}$) are $\tau$-lower
semi-continuous provided card$(\mbox{\rm supp} (\P_{_X}))<\infty$.

For general $X$ we will show that for every  $c \geq 0, \{ D^X_{n,r}> c \}$ is
$\tau$-open. First note that from~(1.4) and~(2.1), there exists a sequence of
quantizations
$\widehat X_m : \Omega \rightarrow F$,
card$(\widehat X_m(\Omega))\le m$,  such that
\[
\lim_m  \| X - \widehat X_m \|_{L_E^r(\P)}=0 .
\]

Consider first the case $r \geq 1$. It follows from
Minkowski's inequality that, for every
$a \in E^n$, $X\mapsto (D^X_{n,r}(a))^{1/r}$ is $1$-Lipschitz on $L_E^r(\P)$:
\begin{eqnarray}
\nonumber | D^X_{n,r}(a)^{1/r} -  D^Y_{n,r} (a)^{1/r} | & = & \left| \|
d(X,a)\|_{L^r(\P)} - \| d(Y, a)\|_{L^r(\P)}\right|\\
\nonumber& \leq & \| d (X, a) - d (Y, a)\|_{L^r(\P)} \\
\label{Lipschitz}  & \leq &  \| X - Y \|_{L_E^r(\P)}.
\end{eqnarray}
Let $a \in \{ (D^X_{n,r})^{1/r} > c\}$. It follows from~(\ref{Lipschitz})
that,  the $\tau$-open set $\{ (D^{\widehat
X_m}_{n,r})^{1/r} > c+ \|X-\widehat X_m\|_{L_E^r(\P)}\}$ is always contained
in $\{ (D^X_{n,r})^{1/r} > c\}$. Furthermore, it contains $a$ for large enough
$m$, still
by~(\ref{Lipschitz}). Hence $\{ (D^X_{n,r})^{1/r}> c\}$ is $\tau$-open and $
D^X_{n,r}$ is $\tau$-l.s.c.

When $0 < r < 1$, one concludes the same way round, using now that
$|u^r-v^r|\le |u-v|^r$ for every $u,\, v\!\in\R_+$, one derives that for every
$a\in E^n$,
\[
| D^X_{n,r}(a) - D^Y_{n,r}  (a) |  \le \E\, | d(X,a) - d(Y, a) |^r\le \| X - Y
\|^r_{L_E^r(\P)}.
\]
\hspace*{\fill}{$\Box$}

In the non-quantization setting $n=1$, Theorem 1 with $F=E$ is due to
Herrndorf (see \cite{HERRN}).

One easily checks that if $E$ is a 1-complemented closed subspace of some
Banach space $G$ and
$(E,G)$ is admissible, then $E$ is admissible. Here $E$ is said to be
$c$-complemented in $G (c \geq 1)$
if there is a linear projection $S$ from $G$ onto $E$ with $\| S \| \leq c$.
An interesting case is $G = E^{**}$. One simply notes that the closed balls in
the bidual $E^{**}$ of $E$
are $\mbox{weak}^*$-compact and thus $E^{**}$ is admissible. The following
characterization is a slight generalization of Theorem 5.9 in \cite{LINDE}.
\begin{Pro}
$(F,E)$ is admissible if and only if
\[
\bigcap_{x \in F} B_E (x, \| z-x \| ) \not= \emptyset \; \mbox{for every} \; z
\in E^{**} .
\]
In particular, if $E$ is 1-complemented in its bidual $E^{**}$, then $E$ is
admissible. \\
\end{Pro}

An investigation of the admissibility feature of Banach spaces $E$ and the
ball topology $\tau_0$ (with
$F = E$) used in the proof of Theorem 1 can be found in \cite{GODEFR},
\cite{GODEFR1}.

One derives for three main classes of Banach spaces the following corollary.
\medskip
\begin{Cor}
In any of the following cases $E$ is 1-complemented in $E^{**}$ and hence, for
every $n \in \N, \; {\cal C}_{n,r}(X,E) \not= \emptyset$.

(i) $E$ is a $KB$ (Kantorovich-Banach)-space.

(ii)  $E$ is a dual space.

(iii) $E$ is an order complete AM-space with unit.
\end{Cor}
{\bf Proof.} $(i)$ By definition, a Banach lattice that is a band in its
bidual is a $KB$-space.
Since $E^{**}$ is an order complete Banach lattice, $E$ is a projection band
in $E^{**}$ and
the band projection from $E^{**}$ onto $E$ has norm 1. (cf. \cite{SCHAE},
Chap. II.5).

$(ii)$ Dual spaces are clearly 1-complemented in their bidual.

$(iii)$ See \cite{SCHAE}, Chap. II.7. \hfill{$\Box$}

\bigskip
The order complete AM-space without unit $c_0 ( \N)$ and the AM-space with
unit $C([0,1])$
which is not order complete admit random variables $X$ without optimal
$n$-quantizers even for $n=1$ (see
the subsequent counterexamples) In particular, both spaces are not admissible.
\\ \\
{\sc Example.} $L^p_{\R^{d}}$-spaces are equipped with the norm
$\| f \|_p = (\int | f (t) | ^p_p d \mu (t) )^{1/p}$ if
$p \in [1, \infty)$ and $\| f \|_{_\infty} = \mu$-ess sup $| f(t) |_{\infty}$
if
$p=\infty$, where $| \cdot |_p$ denotes the $\ell^p$-norm on $\R^d$.
$L^1_{\R^{d}}$-spaces with respect to arbitrary measure spaces and 
are $AL$-spaces and hence $KB$-spaces. $L^p_{\R^{d}}$-spaces, $1 < p < \infty$,
with respect to arbitrary measure spaces 
are reflective and hence dual spaces. $L^\infty_{\R^{d}}$-spaces
with respect to $\sigma$-finite
measure spaces are dual spaces and also order complete
$AM$-spaces with unit
(cf. \cite{SCHAE}, Chap. IV 7). \\ \\
{\bf Remarks.} $\bullet$
Concerning the Banach spaces $E=L^p_{\R^d}$, the above theorem provides new
existence results for the $L^r$-optimal quantizers in the cases $p=1$ and
$p=\infty$. \\
$\bullet$ Any pathwise continuous  process $(X_t)_{t\in[0,1]}$ is an
$L^\infty([0,1], dt)$-Radon  random variable  since $({\cal
C}([0,1]),\|\,.\,\|_{\infty})$ is a Polish subspace of $E=L^\infty([0,1], dt)$
(any probability on a Polish space  is tight $i.e.$ Radon). The above
existence theorem shows that if $\|X\|_{_\infty} \!\in L^r(\P)$ for some
$r>0$, then, for every  $n\ge 1$, $X$ has at least one $L^r$-optimal
$n$-quantizer for the $\|\,.\,\|_{_\infty}$-norm. However, nothing is known
about the pathwise regularity of these optimal quantizers. Surprisingly, we
will see in Section 3 that, for the same process,
$(L^r,\|\,.\,\|_{_p})$-optimal $n$-quantizers with $p < \infty$ have much more
regular paths ($i.e.$ considering $E= L^p$ and $r\ge p$). \\ \\
{\sc Optimal $1$-quantizers may not exist in $c_0(\N)$}
Let $(E,\|\,.\,\|) = (c_0(\N), \|\,.\,\|_{_{\infty}})$ where $c_0(\N)$ denotes
the set of real valued sequences $x =(x_k)_{k\ge 1}$ such that $\lim_k x_k =0$
and $\| x \|_{_{\infty}}= \sup_k | x_k|$.   Let
$(u^{(n)})_{n \geq 1 }$ denote  the canonical basis of
$c_0(\N)$ defined by $u^{(n)}_k = \delta_{n,k}$ where $\delta_{i,j}$ is for
the Kronecker symbol. One considers an $E$-valued  random vector $X$
supported by $\{u^{(n)},\; n\ge 1\}$ with a distribution  $p_n
=\P(X=u^{(n)}),\, n\ge 1$ satisfying $p_n\!\in(0,1/2)$ for every $n\ge 1$. Now
$E^{*} = l^1(\N)$ so that $E^{**} =\ell^\infty(\N)$.
One checks that the assumption of Theorem 1 is not fulfilled either since  the
system
$\{B(u^{(n)},1/2),\, n\ge 1 \}$ has an empty intersection whereas any finite
subsystem has a nonempty intersection.

So  let $n = 1$ and
$r = 1$.  We will show that
$$
e_{1,1} (X, c_0(\N))= 1/2\qquad \mbox{ and }\qquad  {\cal C}_{1,1} (X,
c_0(\N)) = \emptyset.
$$

More precisely we will show that the corresponding level $1$ quantization
problem extended to the Banach space $\ell^\infty(\N)$ does have   a unique
solution $a$ in $\ell^\infty(\N)$ given by
$a_k = 1/2,\, k\ge 1$, that is ${\cal C}_{1,1} (X, \ell^\infty(\N))=\{ a \}$
which in turn implies that it admits no solution in
$c_0(\N)$.   In fact,

\[
\E\,\| X - a \|_{_{\infty}} = \sum\limits^{\infty}_{n = 1} p_n \| u^{(n)} - a
\|_{_{\infty}} = 1/2.
\]
For an arbitrary $b \!\in  \ell^\infty(\N)$ one gets the following: if $\|
u^{(n_{0})} - b \|_{_{\infty}}< 1/2$ for some
$n_0 \geq 1$, then, for every $n \not= n_0$,
\[
\quad\;\;\;\; \| u^{(n)} - b \|_{_{\infty}} \geq   \| u^{(n)} - u^{(n_{0})}
\|_{_{\infty}}  - \| u^{(n_{0})} - b \|_{_{\infty}} = 1 - \| u^{(n_{0})} - b
\|_{_{\infty}}.
\]
Hence
\begin{eqnarray}
\E\, \| X - b \|_{_{\infty}}  & = & \sum _{n \ge 0} p_n  \| u^{(n)} - b
\|_{_{\infty}}   \\
\nonumber & \geq & \sum_{n \not= n_{0}} p_n ( 1 - \| u^{(n_{0})} - b
\|_{_{\infty}}  ) + p_{n_{0}} \| u^{(n_{0})} - b
\|_{_{\infty}}  \\
\nonumber & = & 1 - p_{n_{0}} - ( 1 - 2 p_{n_{0}}) \| u^{(n_{0})} - b
\|_{_{\infty}}  \\
\nonumber & > & 1 - p_{n_{0}} - \frac{1}{2} ( 1 - 2 p_{n_{0}} ) \\
\nonumber & = & 1/2 .
\end{eqnarray}
In case $\| u^{(n)} - b  \|_{_{\infty}}  \geq 1/2$ for every $n \geq 1$, one
clearly obtains
\[
\E \,\| X - b  \|_{_{\infty}} = \sum_{n\ge 1} p_n \| u^{(n)} - b
\|_{_{\infty}}  \geq 1/2 .
\]
According to the above reasoning, any $b \in \ell^\infty(\N)$ that achieves
the infimum must satisfy $\| u^{(n)} - b  \|_{_{\infty}}= 1/2$
for every $n \geq 1$ which  clearly implies $b = a$. Finally
\[
e_{1,1}(X,\ell^\infty(\N))=  \E \| X - a \|_{_{\infty}} = 1/2 \;\mbox{ and
}\;\E \|X- b \|_{_{\infty}}>1/2,\;
a \neq b ,\; b \in \ell^\infty(\N).
\]
On the other hand, as a minimizing sequence from $c_0(\N)$ one may choose
$ a^{(m)} = \frac{1}{2} \sum\limits^{m}_{n=1} u^{(n)} , \;m \geq 1$. Then
\[
\E \| X - a^{(m)} \|_{_{\infty}} = \sum_{n \ge 1} p_n \| u^{(n)} - a^{(m)}
\|_{_{\infty}} = \frac{1}{2} \sum\limits^{m}_{n=1} p_n  +
\sum\limits^{\infty}_{n=m+1} p_n \stackrel{m\to +\infty}{\longrightarrow} 1/2.
\]
Consequently,
\[
e_{1,1} (X, c_0 ( \N)) = 1/2
\]
and since $a \notin c_0(\N)$,  it follows that ${\cal C}_{1,1} (X,c_0(\N))$ is
empty. \\

\smallskip This example  is enlightened by the general Theorem~\ref{E-E**}.
This
theorem solves  the correspondence between  the quantization problem in
$E$ and in $E^{**}$.  It shows that the quantization
error does not decrease when $X$ is seen as random vector in the bidual
$E^{**}$ of $E$ and that the set of its optimal $n$-quantizers as an
$E$-valued random vector is made up with those of its optimal $n$-quantizers
as an $E^{**}$-valued random vector {\em that lie in} $E$. In
particular,
${\cal C}_{n,r} (X,E) = \emptyset$ corresponds to the phenomenon that any
optimal $n$-quantizer of ${\cal C}_{n,r} (X,E^{**})$ has at least one
element in  $E^{**}\setminus E$~: this is precisely what happens in the above
example.

\begin{Thm}\label{E-E**} $(a)$ We have for every $n \in \N$,
\[
e_{n,r} (X,E) = e_{n,r} (X,E^{**}) .
\]
In particular,
\[
{\cal C}_{n,r} (X,E) = \{ \alpha \in {\cal C}_{n,r} (X,E^{**}) : \alpha
\subset E \}.
\]
If card$(\mbox{\rm supp} (\P_{_X})) \geq n$, then $e_{1,r} (X,E) > \cdots >
e_{n,r} (X,E)$.

\smallskip
\noindent $(b)$ Assume that $E$ is admissible. Further assume $\mbox{\rm
supp}(\P_{_X}) = E$. Then
\[
{\cal C}_{n,r} (X,E) = {\cal C}_{n,r} (X,E^{**}) .
\]
\end{Thm}

We first need the following equivariance properties contained in the lemma
below.
\begin{Lem}
Let $E_1$ and $E_2$ be Banach spaces and let $X$ be a Radon $E_1$-valued
random vector satisfying
$\E\| X \|^r < \infty$. If $S : E_1 \rightarrow E_2$ is a bounded linear
operator, then
\[
e_{n,r} (S(X),E_2) \leq \| S \| e_{n,r} (X,E_1) .
\]
If $S : E_1 \rightarrow E_2$ is a bijective linear isometry, $c > 0$ and
$u_2\!\in E_2$, then
\[
e_{n,r} (c\, S(X)+u_2,E_2) = c\, e_{n,r} (X,E_1) \; \mbox{ and } \; {\cal
C}_{n,r} (c\, S (X)+u_2,E_2) = c\, S \left({\cal C}_{n,r} (X,E_1)\right)
+u_2.
\]
\end{Lem}

\noindent
{\bf Proof.} Let us prove $e.g.$ the first assertion. For any $\alpha \subset
E_1$ with $1 \leq \; \mbox{card} \; ( \alpha ) \leq
n$,
\begin{eqnarray*}
e_{n,r} (S(X),E_2) & \leq & ( \E \min_{a \in \alpha} \| S (X) ) - S(a) \|^r
)^{1/r} \\
               & \leq & \| S \| ( \E \min_{a \in \alpha} \| X - a \|^r )^{1/r}
\end{eqnarray*}
and thus the assertion. \hfill{$\Box$}

\bigskip
\noindent {\bf Proof of Theorem~\ref{E-E**}.} $(a)$ The inequality
\[
e_{n,r} (X,E) \geq e_{n,r} (X,E^{**})
\]
is obvious. To prove the converse inequality assume first that $\mbox{\rm
supp}(\P_{_X})$ is finite.
Let $\alpha \in {\cal C}_{n,r}(X,E^{**})$ and let $G$
denote the linear subspace of $E^{**}$ spanned by $\mbox{\rm supp} (\P_{_X})
\cup \alpha$. Since $G$ is
finite-dimensional, there  exists by local reflexivity of $E$, for every
$\varepsilon > 0$, a bounded linear operator $S : G \rightarrow E$ satisfying
$\|S \| \leq 1 + \varepsilon$ and
$S (x) = x$ for every $x \in G \cap E$. (cf. \cite{LITZ} Lemma 1.e.6). Using
Lemma 1, one derives
\begin{eqnarray*}
e_{n,r} (X,E)^r & \leq & \E \min_{b \in S(\alpha)} \| X - b \|^r = \E \min_{a
\in \alpha} \|
S(X) - S (a) \|^r \\ & \leq & (1 + \varepsilon)^r e_{n,r} (X,E^{**})^r .
\end{eqnarray*}
Hence
\[
e_{n,r} (X,E) \leq e_{n,r} (X,E^{**}) .
\]
For general $X$ and $\varepsilon > 0$, choose a quantization $\widehat
X_m:\Omega\to E $ of $X$,
card$(\widehat X_m(\Omega)) \le m$,  for sufficiently large  $m$   such that
\begin{eqnarray*}
\| X - \widehat X_m \|_{L_E^r(\P)}^{1\wedge r} \le  \varepsilon.
\end{eqnarray*}
Then,
\[
| (e_{n,r} (X,E))^{r\wedge 1} - (e_{n,r} (\widehat X_m, E))^{r\wedge 1} |
\leq \varepsilon
\]
and
\[
| (e_{n,r} (X,E^{**}))^{r\wedge 1} - (e_{n,r} (\widehat X_m,E^{**}))^{r\wedge
1} | \leq   \| X - \widehat X_m
\|_{L_E^r(\P)}^{1\wedge r} \le  \varepsilon.
\]
Since card$( \mbox{\rm supp} (\P_{_{\widehat X_m}}) ) \leq m < \infty$, we
have $e_{n,r} (\widehat X_m,E) = e_{n,r} (\widehat X_m,E^{**})$. This
yields
\[
| (e_{n,r} (X,E))^{r\wedge 1} - (e_{n,r} (X,E^{**}))^{r\wedge 1} | \leq 2
\varepsilon .
\]
Hence $e_{n,r} (X,E) = e_{n,r} (X,E^{**})$. Furthermore, since ${\cal C}_{n,r}
(X,E^{**}) \not= \emptyset$ by Corollary 1, it follows from
Proposition~\ref{Pro2.1&2}$(a)$ that $(e_{j,r} (X,E^{**}))_{1 \leq j \leq n}$
is strictly decreasing provided card$( \mbox{\rm supp}(\P_{_X}) )\geq n$.

\smallskip
\noindent $(b)$  The inclusion ${\cal C}_{n,r} (X,E) \subset {\cal C}_{n,r}
(X,E^{**})$ follows from $(a)$. To prove the converse inclusion, we may assume
$\mbox{dim}\,E \geq 1$.
Let $\alpha \in {\cal C}_{n,r} (X,E^{**})$.
By Proposition 2, for every $a \in \alpha$ there exists $b_a \in E$ such that
for every $x \in E$,
\[
\| b_a - x \| \leq \| a-x \|  .
\]
Setting $\beta = \{ b_a : a \in \alpha \}$ this implies
$\beta \in {\cal C}_{n,r} (X,E)$ and that the closed set
\[
A := \{ x \!\in E : \min\limits_{b \in \beta} \| x - b \| = \min\limits_{a \in
\alpha} \| x - a \| \}
\]
satisfies $\P_{_X}(A) = 1$. Therefore, $A = E$ and in  particular, $\beta
\subset A$. One obtains $\min_{a \in \alpha}
\| b - a\| = 0$ for every $b \in \beta$ and hence, $\beta \subset \alpha$. By
Proposition~\ref{Pro2.1&2}$(a)$, we have
card$( \alpha ) = \; \mbox{card} ( \beta ) = n$ which yields $\beta = \alpha$.
Hence $\alpha \in {\cal C}_{n,r} (X,E)$. \hfill{$\Box$} \\ \\
{\bf Remark.} It is to be noticed that the situation  $ {\cal C}_{1,r}
(X,E)=\emptyset$ never occurs for  Gaussian (Radon) random vectors
$X$.  In view of Lemma 1,  we may assume without loss of generality that
$X$ is centered.  Let $r > 0$. It  follows from the Anderson inequality
(\cite{LETA}) that, for every $a \!\in E$,
\[
\E \|X- a \|^r = \int_0^{+\infty} \P (\|X- a \|^r\ge t )d t  \ge
\int_0^{+\infty} \P (\|X\|^r\ge t )d t =\E \|X\|^r
\]
so that    $\{ 0 \} \in {\cal C}_{1,r} (X,E)\neq \emptyset$. However, it
remains an open question whether
${\cal C}_{n,r} (X,E)$ may be empty for  $n \geq 2$ or not. \\

An immediate consequence of Theorem 2(a) is as follows. Let us call a Banach
subspace $F$ of $E$ {\em locally $c$-complemented} $(c \geq 1)$ if there is a
linear operator $S : E \rightarrow F^{**}$ of norm $\| S\|  \leq c$
satisfying $S(x) = x$ for every $x \in F$. Notice that local 1-complementation
coincides with the notion of an ideal introduced  in \cite{GODEFR2}.
\begin{Cor}
Assume that $\P_X (F) =1$ for some Banach subspace $F$ of $E$ and that $F$ is
locally 1-complemented in $E$. Then, for every $n \in \N$,
\[
e_{n,r} (X,F) = e_{n,r} (X,E) .
\]
In particular, ${\cal C}_{n,r} (X,F) \not= \emptyset$ implies ${\cal
C}_{n,r}(X,E) \not= \emptyset$. \\ 
\end{Cor}
{\bf Proof.} It follows from Theorem 2(a) and Lemma 1 that
\[
\begin{array}{lcl}
e_{n,r} (X,F) & = & e_{n,r} (X, F^{**}) = e_{n,r} (S(X),F^{**}) \\
& \leq & \| S \| e_{n,r} (X,E) = e_{n,r} (X,E) \leq e_{n,r} (X,F) .
\end{array}
\]
\hspace*{\fill}{$\Box$}

\bigskip
One observes that the preceding corollary contains Theorem 2(a) since $E$ is
obviously locally 1-complemented in $E^{**}$.

\bigskip
\noindent {\sc Example} $\bullet$ AM-spaces $F$ are locally 1-complemented as
Banach subspace in any Banach space $E$. In fact, since $F^{**}$ is an
order complete AM-space with unit, this feature follows from Theorem II.7.10
in \cite{SCHAE}. For instance, if $E = C(T)$ for some compact metric
space $T$ and
\[
F = \{ f \in C(T) : f(t) = 0 \; \mbox{for all} \; t \in T_0\}
\]
for some closed subset $T_0$ of $T$, then $F$ is a closed vector sublattice of
the AM-space $C(T)$ and thus an AM-space.  \\
$\bullet$ AL-spaces $F$ are 1-complemented as Banach sublattice in any Banach
lattice $E$ (see \cite{SCHAE}, II.8). \\

Finite dimensional subspaces of dimension $d \geq 2$ are admissible but not
necessarily (locally) 1-complemented. In fact, it may happen that $\P_X(F) =
1$ for some 2-dimensional subspace $F$ of $E$ and
${\cal C}_{n,r}(X,E) = \emptyset$ even for $n=1$. In particular, $(F,E)$ is
not admissible.
The following example is taken from Herrndorf \cite{HERRN}. \\ \\
{\sc A counterexample when dim $F=2$} Let $\ell^1(\N)$ be equipped with the
$\ell^1$-norm
$\|  x \| = \sum\limits^{\infty}_{j=1} |  x_j | $. Let
$(u^{(n)})_{n \geq 1}$ be the canonical basis of $\ell^1(\N)$ and set
$v^{(1)} : = 0, v^{(2)} := u^{(1)} - u^{(2)}$ and $v^{(3)} := u^{(1)} -
u^{(3)}$. Consider the $\ell^1(\N)$-valued random variable $X$ supported by
$\{v^{(1)}, v^{(2)}, v^{(3)} \}$ with $\P(X = v^{(i)}) = 1/3$. Let $F$ denote
the linear span of $\{ v^{(2)}, v^{(3)} \}$ in $ \ell^1 (\N)$. So $\P (X \in
F) = 1$ and dim$F = 2$.

Let $n = 1$ and $r = 1$. First will show that
\[
e_{1,1} (X,F) = 4/3 , \; e_{1,1} (X, \ell^1 ( \N )) = 1
\]
and
\[
{\cal C}_{1,1} (X, \ell^1 ( \N)) = \{ \{ u^{(1)} \} \}.
\]
In fact,
\[
 \E \| X - u^{(1)} \| = \frac{1}{3} \sum\limits^{3}_{i=1} \| v^{(i)} - u^{(1)}
\| = 1 .
\]
On the other hand, once noticed that $\| v^{(i)} - v^{(j)} \|  = 2$ for $i
\not= j$,
one shows like in the previous counterexample that for every
$a \in \ell^1 (\N), \E \| X - a\| = \frac{1}{3} \sum\limits^{3}_{i=1} \|
v^{(i)} - a \|  \geq 1$ and that any $L^1$-optimal 1-quantizer 
$a \in \ell^1 (\N)$ must satisfy
$\| v^{(i)} - a \| = 1$ for every $i \in \{1, 2, 3 \}$ which implies $a =
u^{(1)}$.
As for $e_{1,1} (X,F)$, observe that
\[
\E \|  X - v^{(i)} \| = 4/3 , \; i \in \{ 1, 2, 3 \} .
\]
Any $a \in F$ can be written as $a = (s+t) u^{(1)} - s u^{(2)} - t u^{(3)}, s,
t \in \R$, so that
\[
\begin{array}{l c l}
\sum\limits^{3}_{i=1} \| v^{(i)} - a \|
& = & | s+t|  + | s|  + | t| + | 1-s-t|  + | 1-s|  + | t| \\
&   & + | 1-s-t|  + | s|  + | 1-t|  \\
& \geq & 4
\end{array}
\]
since $| 1-t|  + | t|  \geq 1, t \in \R$. This yields $e_{1,1} (X, F) = 4/3$.

Now we construct a Banach subspace $E$ of $\ell^1 (\N)$ such that $F \subset
E$ and
\[
{\cal C}_{1,1} (X, E) = \emptyset .
\]
Choose $c = (c_j)_{j \geq 1} \in \ell^\infty ( \N)$ such that $c_1 = c_2 = c_3
= 1$ and
$(c_j)_{j \geq 3}$ is strictly increasing with
$\| c \|_{\infty} = \sup_{j \geq 1} | c_j |  > 3$. Define $E$ as the
hyperplane
\[
E := \{x \in \ell^1 (\N) : \sum\limits^{\infty}_{j=1} x_j c_j = 0 \} .
\]
Then $F \subset E$. For $k \geq 4$, set $a^{(k)} := u^{(1)} - \frac{1}{c_{k}}
u^{(k)} $. One obtains
$a^{(k)} \in E$ and
\[
\E \| X - a^{(k)} \|  = \frac{1}{3} \sum\limits^{3}_{i=1} \| v^{(i)} - a^{(k)}
\|
= 1 + 1/c_k .
\]
Consequently,
\[
e_{1,1} (X,E) \leq 1 + 1/\| c \|_{\infty} < 4/3 = e_{1,1} (X,F) .
\]
For an arbitray $a \in E$ one gets the following: if $a_j = 0$ for $j \geq 4$,
then $a \in F$
and hence $\E \| X - a \| \geq 4/3 > e_{1,1} (X,E)$. If
$a_j \not= 0$ for some $j \geq 4$, a can be strictly improved. Set
\[
b := a - a_j u^{(j)} + a_j c_j c^{-1}_{j+1} u^{(j+1)} .
\]
One checks that $b \in E$ and for every $i \in \{1, 2, 3\}$,
\[
\| v^{(i)} - b \|  = \sum\limits^{3}_{k=1} | v^{(i)}_k - a_k | +
\sum\limits_{ { {k \geq 4}\atop{k \not= j, j + 1}} } | a_k | + | b_j | + |
b_{j+1} |
< \|  v^{(i)} - a \|  .
\]
This implies
\[
\E \| X - b \|  < \E \| X - a \|  .
\]
Consequently, ${\cal C}_{1,1} (X, E) = \emptyset$.
\subsection{Optimal quantizers for continuous stochastic processes}
Now we turn to $\R^d$-valued pathwise continuous processes $X = (X_t)_{t \in
T}$
indexed by a compact metric space $T$. The space $E := C_{\R^{d}}(T)$ of
$\R^d$-valued continuous
functions on $T$ and the space $M^b_{\R^{d}}(T)$ of bounded, $\R^d$-valued,
Borel measurable
functions on $T$ are Banach spaces under the norm
\begin{equation}
\| f \|_{\sup} := \sup_{t \in T}| f(t)|_{_\infty}
\end{equation}
where $| \cdot |_{_\infty}$ denotes the $\ell^\infty$-norm on $\R^d$. Since
$C_{\R^{d}}(T)$
is separable, $X$ is Radon when viewed as $C_{\R^{d}}(T)$-valued random
variable. Consequently,
$X$ is Radon as $M^b_{\R^{d}}(T)$-random variable.
\begin{Thm} Let $T$ be compact metric space. Then the pair
$(C_{\R^{d}}(T), M^b_{\R^{d}}(T))$ is admissible under the norm (2.9). In
particular, if
$X = (X_t)_{t \in T}$ is a $\R^d$-valued pathwise continuous process with
$\E \| X \|^r_{sup} < \infty$, then for every $n \in \N$,
\end{Thm}
\[
{\cal C}_{n,r} (X, M^b_{\R^{d}}(T)) \not= \emptyset .
\]

The proof of Theorem 3 is based on the admissibility of
$L^\infty_{\R^{d}}$-spaces and the
following ``lifting property''.

\begin{Lem} Let $\mu$ be a finite Borel measurable on the compact metric space
$T$ with
$\mbox{supp} (\mu) = T$. Then for every $h \in M^b_{\R^{d}}(T)$ there exists
$g \in M^b_{\R^{d}}(T)$ such that $g = h \; \mu$-a.e. and
\[
\| f-g \|_{\sup} = \| f-h \|_{_\infty} \; \mbox{for every} \;
f \in C_{\R^d}(T)
\]
where
\end{Lem}
\begin{equation}
\| h \|_{_\infty} := \mu\mbox{-esssup} \; | h |_{_\infty} .
\end{equation}
{\bf Proof.} One notes that for $h = (h_1 , \ldots , h_d) \in M^b_{\R^{d}}
(T)$,
\[
\| h \|_{\sup} = \max_{1 \leq i \leq d} \| h_i \|_{\sup}
\]
and
\[
\| h \|_{_\infty} = \max_{1 \leq i \leq d} \| h_i \|_{_\infty} .
\]
Therefore, it is enough to consider the case $d = 1$. Set $C(T) = C_\R(T)$ and
$M^b(T) = M^b_\R(T)$. Let $D$ be a countable dense subset of $C(T)$. Observe
that the norms
$\| \cdot \|_{_\infty}$ and $\| \cdot \|_{\sup}$ coincide
on $C(T)$. This is a consequence of the
assumption $\mbox{supp}(\mu) = T$. Let $h \in M^b(T)$. For $f \in C(T)$, set
\[
c_f := \| f-h \|_{_\infty} .
\]
Then
\[
N^{+}_f := \{t \in T : f(t) - h(t) > c_f\}
\]
and
\[
N^{-}_f := \{t \in T : h(t) - f(t) > c_f\}
\]
are Borel subsets of $T$ with $\mu$-measure zero. Consequently,
\[
N := \bigcup\limits_{f \in D} ( N^{+} _{f} \cup N^{-}_f)
\]
satisfies $\mu(N) = 0$. Since for every $t \in T \setminus N$ and $f \in D$,
\[
f(t) - h(t) \leq c_f \; \mbox{ and } \; h(t) - f(t) \leq c_f
\]
one obtains
\begin{equation}
\sup_{t \in K \setminus N} | f(t) - h(t) | \leq c_f , f \in D .
\end{equation}
The construction of the function $g$ is given in two steps. \\ \\
{\sc Step 1.} For $\varepsilon > 0$ and $t \in T$, let
\[
d(t, \varepsilon) := \mu \mbox{-esssup} \; h_{| U (t, \varepsilon)} ,
\]
where $U(t, \varepsilon)$ denotes
the open ball in $T$ of radius $\varepsilon$ centered at $t$.
Define the ``upper limit function'' $\hat{h} : T \rightarrow \R$ of $h$ by
\[
\hat{h} (t) := \lim\limits_{\varepsilon \downarrow 0} d(t, \varepsilon) .
\]
One easily checks that for any Borel subset $A$ of $T$, the function
$T \rightarrow \R, t \mapsto \mu (U(t, \varepsilon) \cap A)$ is Borel.
Therefore, for every $a \in \R$,
\begin{eqnarray*}
    \{ t \in T : \hat{h} (t) < a \} & =&\left\{ t \in T : \exists\, n \in \N,
\exists\, m \in \N \; \mbox{such that} \;
       h | U( t, \frac{1}{n}) \leq a - \frac{1}{m} \; \mu\mbox{-a.e.} \right\}
\\
 & = & \bigcup_{n \in \N} \bigcup_{m \in \N} \left\{ t \in T : \mu ( U ( t,
      \frac{1}{n}) \cap \{h > a - \frac{1}{m} \} ) = 0 \right\}
\end{eqnarray*}
is a Borel set and thus $\hat{h}$ is Borel measurable. The function $\hat{h}$
has the following property: for every $t \in N$ there exists a sequence
$(t_n)$ in $T \setminus N$ such that $\lim\limits_{n \to \infty} t_n = t$ and
$\lim_{n \to \infty} h(t_n) = \hat{h} (t)$. In fact, let $t \in N$ and let
$\varepsilon_n \downarrow 0$
so that $\hat{h}(t) = \lim_{n \to \infty} d(t, \varepsilon_n)$.
For every $n \in \N$, there exists $t_n \in U(t, \varepsilon_n) \setminus N$
such that
\[
d(t, \varepsilon_n) - \frac{1}{n} < h(t_n) \leq d(t, \varepsilon_n).
\]
This implies
\[
\lim\limits_{n \to \infty} h(t_n) = \hat{h}(t) \; \mbox{and} \; \lim\limits_{n
\to \infty} t_n = t .
\]
{\sc Step 2.} Define $g : T \rightarrow \R$ by
\[
g(t) := \left\{ \begin{array}{ll}
\hat{h} (t) & , t \in N \\
h(t) & , t \in T \setminus N .
\end{array} \right.
\]
We show that $g$ has the required properties. Observe that $g$ is Borel
measurable,
$g = h \; \mu$-a.e. and $\| g \|_{\sup} \leq \| h \|_{sup} < \infty$.
Let $f \in D$   If $t \in T \setminus N$, then $g (t) = h(t)$ and hence by
(2.4),
$| f(t) - g(t) | \leq c_f$. By step 1, if $t \in N$, there exists a sequence
$(t_n)$ in
$T \setminus N$ such that $\lim t_n = t$ and $\lim h(t_n) = \hat{h} (t)$.
Therefore,
\begin{eqnarray*}
| f (t) - g(t)|  & = & | f(t) - \hat{h} (t)| = | \lim\limits_{n \to \infty}
f(t_n) - \lim_{n \to \infty} h(t_n) | \\
 & = & \lim\limits_{n \to \infty} | f(t_n) - h(t_n)|  \\
 & \leq & \sup\limits_{s \in T \setminus N} | f(s) - h(s) | \leq c_f .
\end{eqnarray*}
Consequently,
\begin{equation}
\| f-g \|_{\sup} \leq c_f, f \in D .
\end{equation}
Now let $f \in C(T)$. There exists a sequence $(f_n)$ in $D$ such that
$\lim\limits_{n \to \infty} \| f-f_n \|_{\sup} = 0$. For every $ t \in T$,
\begin{eqnarray*}
| f(t) - g(t) | & = & | \lim\limits_{n \to \infty} f_n (t) - g(t) | =
\lim\limits_{n \to \infty} | f_n(t) - g(t) | \\
 & \leq & \limsup_{n \to \infty} c_{f_{n}} .
\end{eqnarray*}
Since
\[
c_{f_{n}} = \| f_n - h \|_{_\infty} \leq \| f_n - f \|_{\sup} + c_f
\]
one obtains
\begin{equation}
\| f-g \|_{\sup} \leq \limsup_{n \to \infty} c_{f_{n}} \leq c_f .
\end{equation}
Conversely, we clearly have
\[
c_f = \| f-h \|_{_\infty} = \| f-g \|_{_\infty} \leq \| f-g
\|_{\sup}.
\]
\hspace*{\fill}{$\Box$}

\bigskip
\noindent  {\bf Proof of Theorem 3.} Let ${\cal K} = \{B_{M^b} (f_i, \rho_i) : i \in
I\}$ be a system of closed
balls in $M^b_{\R^{d}} (T)$ with centers $f_i \in C_{\R^{d}} (T)$ satisfying
the finite
intersection property. Choose a finite Borel measure $\mu$ on $T$ such that
$\mbox{supp}(\mu)= T$
and consider the system $\tilde{\cal K} = \{ B_{L^\infty} (S f_i, \rho_i) : i \in
I\}$ of corresponding
closed balls in $L^\infty_{\R^{d}}(\mu)$ under the norm $\|  \cdot
\|_{\infty}$
(see (2.10)) where 
$S : M^b_{\R^d} (T) \rightarrow L^\infty_{\R^d} (\mu)$ denotes the quotient map. 
It is obvious that $\tilde{\cal K}$ also has the finite
intersection property. Since
$L^\infty_{\R^{d}}  (\mu)$ is admissible by Proposition 2 and
Corollary 1, $\tilde{\cal K}$ has a nonempty
intersection. Let $S(h)$ be a member of this
intersection. Lemma 2 implies that there is a function
$g \in M^b_{\R^{d}} (T)$ such that $g = h \; \mu$-a.e. and
\[
\| f_i - g \|_{\sup} = \| f_i - h \|_{_\infty} = \parallel S f_i - S h \parallel_\infty \; \mbox{for every} \;
i \in I .
\]
Consequently, $g$ belongs to the intersection of ${\cal K}$. This yields the
required admissibility.
\hfill{$\Box$}

\bigskip

One derives from Corollary 2 that the quantization error
does not decrease when $X$ is seen as $M^b_{\R^d} (T)$-or even
$L^\infty_{\R^d}(\mu)$-valued random variable. \\ \\
{\bf Theorem 4} {\em Assume that $X = (X_t)_{t \in T}$ is a $\R^d$-valued
pathwise continuous
process indexed by a compact metric space $T$ with $\E \| X \|^r_{sup} <
\infty$. .
Let $\mu$ be a finite Borel measure on $T$ with
$\mbox{supp} ( \mu ) = T$. Then for every $n \in \N$,
\[
e_{n,r} (X, C_{\R^{d}}(T)) = e_{n,r} (X, M^b_{\R^{d}} (T)) = e_{n,r}
(X, L^\infty_{\R^{d}} ( \mu)) ,
\]
where $C_{\R^{d}} (T)$ and $M^b_{\R^{d}}(T)$ are equipped with the sup-norm
(2.9) and
$L^\infty_{\R^{d}}(\mu)$ is equipped with the norm (2.10). In particular, }
\begin{eqnarray*}
{\cal C}_{n,r} (X,C_{\R^{d}}(T)) & = & \{\alpha \in {\cal C}_{n,r} (X,
M^b_{\R^{d}}(T)) : \alpha
\subset C_{\R^{d}} (T) \}  \\
& = & \{ \alpha \in {\cal C}_{n,r} (X, L^\infty_{\R^d}(\mu)) : (\mbox{a} \;
\mu-\mbox{version of} \;) \alpha \subset C_{\R^d}(T) \} .
\end{eqnarray*}
{\bf Proof.}
$C_{\R^d} (T)$ is an AM-space so that Corollary 2 applies. We obtain
\[
e_{n,r} (X, C_{\R^{d}}(T)) = e_{n,r} (X, M^b_{\R^{d}}(T)) .
\]
Since $C_{\R^d} (T)$ can be considered as a subspace of $L^{\infty}_{\R^d} (\mu)$, the same argument yields
\[
e_{n,r} (X, C_{\R^{d}} (T)) = e_{n,r} (X, L^\infty_{\R^{d}}(\mu)) .
\]
(The latter equality is also an immediate consequence of Lemma 2.) \hfill{$\Box$} \\

We will exhibit a pathwise continuous process $X = (X_t)_{t \in
[0,1]}$ having no $L^1$-optimal
1-quantizer in $C([0,1])$. In particular, due to the lack of order
completeness, $C([0,1])$ is not
admissible. \\ \\
{\sc Optimal 1-quantizer may not exist in $C([0,1])$}
Let $(E, \| \cdot \| )=(C([0,1]), \| \cdot \|_{\sup})$.
Define, for every $n \in \N$, a continuous function $f_n : [0,1] \rightarrow
\R$ by
\[
f_n(t) := \left\{
\begin{array}{lcl}
0  & \mbox{if} & t \in [0, \frac{1}{2} - 2^{-n}] \cup [\frac{1}{2} -
2^{-(n+1)} , \frac{1}{2} ] \\
2^{n+1} (2t-1)+4 & \mbox{if} & t \in [ \frac{1}{2} - 2^{-n} , \frac{1}{2} - 3
\cdot 2^{-(n+2)}] \\
2^{n+1} (1-2t) - 2 & \mbox{if} & t \in [ \frac{1}{2} - 3 \cdot 2^{-(n+2)},
\frac{1}{2} - 2^{-(n+1)}] \\
- f_n (1-t) & \mbox{if} & t \in [\frac{1}{2} , 1 ] .
\end{array} \right.
\]
One considers an E-valued random variable $X$ supported by $\{f_n : n \geq 1
\}$ with
$p_n := \P (X = f_n)$ satisfying $p_n \in (0, 1/2)$ for every $n \in \N$ and
$\sum\limits^{\infty}_{n=1} p_n = 1$. The assumption of Theorem 1 is not
fulfilled since the
system $\{B_E (f_n, \frac{1}{2}) : n \geq 1 \}$ has the finite intersection
property whereas it
has an empty intersection (see below).

Let $n = 1$ and $r = 1$. We will show that
\[
e_{1,1} (X, E) = 1/2 \; \mbox{and} \; {\cal C}_{1,1} (X, E) = \emptyset .
\]
Recall that by Theorem 4, $e_{1,1} (X,E) = e_{1,1} (X,G)$ where $G =
M^b([0,1])$ equipped with
$\| \cdot \|_{\sup}$. Set $h := \frac{1}{2} ( 1_{[0, 1/2]} -1_{(1/2, 1]} )$.
One checks that, for every $n \geq 1$,
\[
\| f_n - h \|_{\sup} = 1/2
\]
so that
\[
\E \| X-h \|_{\sup} = \sum^\infty_{n=1} p_n \| f_n - h \|_{\sup} = 1/2.
\]
On the other hand, one shows like in the $c_0 (\N)$-counterexample preceding
Theorem 2 that for every
$g \in G, \E \| X-g\|_{\sup} \geq 1/2$ and that any $L^1$-optimal 1-quantizer
$\{ g \}$
must satisfy $\| f_n - g \|_{\sup} = 1/2$ for every $n \in \N$: one reproduces
the string of
inequalities starting at (2.8) once noticed that $\| f_n - f_m \|_{\sup} = 1$
for every
$n \not= m$.
This implies $e_{1,1} (X,G) = 1/2$ and $\{ h \} \in {\cal C}_{1,1} (X,G)$.
Furthermore, no $g \in E$ can satisfy the condition $\| f_n - g \|_{\sup} =
1/2$ for every
$n \in \N$.
In fact, if $g (1/2) < 1/2$, then $g(t_n) < 1/2$ with
$t_n = \frac{1}{2} - 3 \cdot 2^{-(n+2)}$ and $n$ large enough so that
\[
| f_n (t_n) - g (t_n) | = 1 - g(t_n) > 1/2 .
\]
If $g(1/2) \geq 1/2$, then $g (1 - t_n) \geq 0$ for $n$ large enough so that
\[
| g (1 - t_n) - f_n (1 - t_n) | = g (1-t_n) + 1 \geq 1 .
\]
Consequently, ${\cal C}_{1,1} (X,E) = \emptyset$.
\subsection{Bounds for quantization errors}
As before let $X$ be a Radon random variable in $(E, \| \cdot \| )$ satisfying
the integrability condition (1.2). The following observation (a) is already
contained in \cite{CREUT}. 
\begin{Pro}
Assume that $\P_X(F) = 1$ for some Banach subspace $F$ of $E$. \\
(a) For every $n \in \N$,
\[
e_{n,r} (X,E) \leq e_{n,r} (X,F) \leq 2 e_{n,r} (X,E) .
\]
(b) If $F$ is locally c-complemented in $E$, then for every $n \in \N$,
\[
e_{n,r} (X,F) \leq c e_{n,r} (X,E).
\]
\end{Pro}
{\bf Proof.} (a) We have to prove only the second inequality. Let
$\alpha = \{a_1, \ldots, a_n \} \subset E$ and $\varepsilon > 0$. Choose $b_i
\in F$ such that
$\| a_i - b_i \| \leq ( 1 + \varepsilon) \mbox{dist} (a_i , F)$.
This implies that
\[
\| a_i - b_i \| \leq (1 + \varepsilon) \| X - a_i \| \; \mbox{a.e.}
\]
for every $i \in \{ 1, \ldots , n\}$ and hence
\[
\min_{1 \leq i \leq n} \| X - b_i \| \leq (2 + \varepsilon) \min_{1 \leq i
\leq n}
\| X - a_i \| \; \mbox{a.e.}
\]
Consequently,
\[
e_{n,r} (X,F) \leq (2 + \varepsilon) (\E \min_{a \in \alpha} \| X-a
\|^r)^{1/r} .
\]
This yields the assertion.

\smallskip
\noindent (b) is an immediate consequence of Theorem 2(a) and Lemma 1.
\hfill{$\Box$} \\

It is to be noticed that the factor 2 in part (a) of the preceding proposition
is sharp. It cannot be improved as universal constant.
This is demonstrated in the subsequent Example. In view of (a), the cases of
interest in part (b) are
$c < 2$. \\ \\
{\sc The constant 2 is sharp}
We modify the setting of the counterexample following Corollary 2. Let
$E = \ell^1 ( \N ), \| x \|  = \sum\limits^{\infty}_{j=1} |  x_j | $
and let $(u^{(n)})_{n \geq 1}$ denote the canonical basis of E. Fix $m \in \N,
m \geq 2$ and set
$v^{(i)} := u^{(1)} - u^{(i)}, i \in \{ 1, \ldots , m \}$. One considers the
$E$-valued random variable $X$
supported by $\{v^{(1)} , \ldots , v^{(m)} \}$ with $\P (X = v^{(i)})= 1/m$.
Let $F$ denote the linear
span of $\{ v^{(1)}, \ldots , v^{(m)} \}$. So $\P (X \in F) = 1$.

Let $n = 1$ and $r =1$. One checks like in the above mentioned counterexample
that
\[
e_{1,1} (X, E) = 1 .
\]
We will show that
\[
e_{1,1} (X, F) = 2 (m-1)/m .
\]
In fact, for $j \in \{1 , \ldots , m \}$,
\[
\E \| X - v^{(j)} \| = \frac{1}{m} \sum\limits^{m}_{i=1} \| v^{(i)} - v^{(j)}
\|
= 2 (m-1)/m .
\]
Any $a \in F$ can be written as
$a = \sum\limits^{m}_{j=2} s_j u^{(1)} - \sum\limits^{m}_{j=2} s_j u^{(j)},
s_j \in \R$
and hence
\[
\begin{array}{lcl}
\| v^{(1)} - a \|  & = & \|  a \| =  |  \sum\limits^{m}_{j=2} s_j |
+ \sum\limits^{m}_{j=2} | s_j |  , \\
\| v^{(i)} - a \|  & = & | 1 - \sum\limits^{m}_{j=2} s_j | + |  1 - s_i |
+ \sum\limits^{m}_{\stackrel{j=2}{j\not= i}} | s_j | , i \in  \{ 2, \ldots , m
\} .
\end{array}
\]
Using the elementary inequalities $| 1-t |  + | t |  \geq 1$ and
$| 1-s-t |  + | s | + | t |  \geq 1, s, t \in \R$, one obtains
\[
\sum\limits^{m}_{i=1} \| v^{(i)} - a \|  \geq 2(m-1) .
\]
Consequently,
\[
\E \| X - a \|  \geq 2(m-1)/m . \vspace{0.5cm}
\]

Next we describe marginal bounds for $\R^d$-valued stochastic processes. For
$p \in [1, \infty)$, let $E = L^p_{\R^d}(T, {\cal B}, \mu)$,
$\mu$ finite measure,  equipped with the norm
\begin{equation}
\| f \|_p := ( \int | f(t) |^p_p d \mu(t))^{1/p} = ( \sum\limits^d_{i=1} \int
| f_i (t) |^p d\mu(t) )^{1/p} .
\end{equation}
Assume that $E$ is separable.
Let $X = (X_t)_{t \in T} = (X_{1,t} , \ldots X_{d,t})_{t \in T}$ be a
bi-measurable $\R^d$-valued process such that
\begin{equation}
\E \| X \|^p_p < \infty .
\end{equation}
Then the process $X$ can be seen as a (Radon) random vector taking its values
in $L^p_{\R^d}$. For the sake of
simplicity, we consider the case $r=p$. As for bounds when constants are not
important there will be no loss of generality since usual inequalities on
$L^p$-norms imply for $r \in [1, \infty)$
\[
\mu(T)^{ \frac{1}{p}- \frac{1}{p \wedge r}} e_{n,p \wedge r} (X, L^{p \wedge
r}_{\R^d} ) \leq e_{n,r}
(X, L^p_{\R^d}) \leq \mu(T)^{\frac{1}{p}-\frac{1}{p \wedge r}} e_{n, p \vee r}
(X, L^{p \vee r}_{\R^d}).
\]
\begin{Pro}
Let $p \in [1, \infty)$. For every $n, n_1, \ldots , n_d \in \N$ such that
$\prod\limits^{d}_{i=1} n_i \leq n$,
\[
\sum\limits^{d}_{i=1} e_{n,p} (X_i, L^p)^p \leq e_{n,p} (X,L^p_{\R^d})^p \leq
\sum\limits^{d}_{i=1} e_{n_i, p} (X_i , L^p)^p . \vspace{1.0cm}
\]
\end{Pro}
{\bf Proof.} As for the upper estimate, let $\alpha_i \subset L^p$ be a
$L^p$-optimal $n_i$-quantizer for
$X_i, i \in \{1 , \ldots , d\}$ (see Corollary 1). Set $\alpha :=
\times^{d}_{i=1} \alpha_i$. Thus $\alpha$ consists of functions
$a = (a_1 , \ldots , a_d) \in L^p_{\R^d}$ with $a_i \in \alpha_i$ and
card$(\alpha) \leq n$. One obtains
\[
\begin{array}{lcl}
e_{n,p} (X, L^p_{\R^d})^p & \leq & \E \min_{a \in \alpha} \| X-a \|^p_p \\
& = & \E \min_{a \in \alpha} \sum\limits^{d}_{i=1} \int | X_{i,t} - a_i(t) |^p
d\mu(t) \\
& = & \E \sum\limits^{d}_{i=1} \min_{b \in \alpha_i} \int | X_{i,t} - b(t) |^p
d\mu(t) \\
& = & \sum\limits^{d}_{i=1} e_{n_i,p} (X_i, L^p)^p .
\end{array}
\]
As for the lower estimate, let $\alpha \subset L^p_{\R^d}$ with card$(\alpha)
\leq n$. Then
\[
\begin{array}{lcl}
\E \min_{a \in \alpha} \| X -a \|^p_p & \geq & \E \sum\limits^d_{i=1} \min_{a
\in \alpha} \int | X_{i,t} - a_i (t)|^p d\mu(t) \\
& \geq & \sum\limits^d_{i=1} e_{n,p}(X_i, L^p)^p .
\end{array}
\]
This yields the lower estimate. \hfill{$\Box$}

\bigskip
Now let $T$ be a compact metric space and assume that $X = (X_t)_{t \in T}$ is
a $\R^d$-valued continuous process. Let $E =C_{\R^d}(T)$
equipped with the sup-norm (2.9). Assume
\begin{equation}
\E \| X \|^r_{\sup} < \infty .
\end{equation}
\begin{Pro}
Let $r \in (0, \infty)$. Let $c \in (0, \infty)$ such that $| \cdot
|_{_\infty} \leq c | \cdot |_r$.
Then for every $n, n_1 , \ldots , n_d \in \N$ such that $\Pi^d_{i=1} n_i \leq
n$,
\[
\max_{1 \leq i \leq d} e_{n,r} (X_{i}, C(T))^r \leq e_{n,r} (X,C_{\R^d} (T))^r
\leq c^r \sum\limits^d_{i=1}
 e_{n_i,r} (X_i, C(T))^r .
\]
\end{Pro}
{\bf Proof.} For $i \in \{ 1, \ldots , d\}$ and $\varepsilon > 0$, choose
$\alpha_i \subset C(T)$ such that card$(\alpha_i) \leq n_i$ and
\[
\E \min_{b \in \alpha_i} \| X_i - b\|^r_{\sup} \leq e_{n_i,r} (X_i , C(T))^r +
\varepsilon .
\]
Set $\alpha := \times^d_{i=1} \alpha_i$. Then $\alpha \subset C_{\R^d}(T)$,
card$(\alpha) \leq n$ and
\[
\begin{array}{lcl}
\E \min_{a \in \alpha} \| X-a \|^r_{\sup} & \leq & c^r \E \min_{a \in \alpha}
\sup_{t \in T}
          \sum\limits^d_{i=1} | X_{i,t} - a_i(t) |^r \\
& \leq & c^r \E \min_{a \in \alpha} \sum\limits^d_{i=1} \| X_i - a_i
\|^r_{\sup} \\
&   =  & c^r \E \sum\limits^d_{i=1} \min_{b \in \alpha_i} \| X_i - b
\|^r_{\sup} \\
& \leq & c^r \sum\limits^d_{i=1} e_{n_i,r} (X_i, C(T))^r + c^r d \varepsilon .
\end{array}
\]
This yields the upper esxtimate. As for the lower estimate, let $\alpha
\subset C_{\R^d} (T)$ with card$(\alpha) \leq n$. Then for every $i$,
\[
\E \min_{a \in \alpha} \| X-a \|^r_{\sup} \geq \E \min_{a \in \alpha} \| X_i -
a_i \|^r_{\sup} \geq e_{n,r} (X_i, C(T))^r
\]
which gives the lower estimate. \hfill{$\Box$} \\

In the preceding proposition one may replace $C_{\R^d}(T)$ and $C(T)$ by
$L^\infty_{\R^d} (\mu)$ and
$L^\infty(\mu)$ respectively for any finite Borel measure $\mu$ on $T$ with
$\mbox{supp} (\mu) = T$. This follows from Theorem 4. \\ \\
\section{Stationary quantizers}
\setcounter{equation}{0}
\setcounter{Assumption}{0}
\setcounter{Theorem}{0}
\setcounter{Proposition}{0}
\setcounter{Corollary}{0}
\setcounter{Lemma}{0}
\setcounter{Definition}{0}
\setcounter{Remark}{0}
Let $X$ be a Radon $(E, \| \cdot \| )$-valued random variable satisfying
condition (1.2).
We will introduce a notion of $L^r$-stationary quantizer as the critical
points of level $n$ $L^r$-distortion function $D_{n,r}^X$
formerly defined by Equation~(2.6).
For a quantizer $\alpha = \{ a_1 , \ldots , a_n \}$ let $V_i (\alpha) =
V_{a_{i}}(\alpha)$ and $C_i(\alpha) = C_{a_{i}}(\alpha)$.
\begin{Dfn}   A $n$-quantizer $\alpha = \{a_1, \dots, a_n\}\subset E$ of size
$n$ is called admissible for
$X$ if
\[
\left\{\begin{array}{ll}
(i) &\P_{_X}(V_i(\alpha))>0,\quad i=1,\ldots,n,\\
*[.4em]
(ii) & \P_{_X}(V_i(\alpha)\cap V_j(\alpha))=0,\quad i,\,j=1,\ldots,n,\; i\neq
j.
\end{array}\right.
\]
A $n$-tuple $(a_1,\ldots,a_n)\!\in E^n$ is admissible if its associated
$n$-quantizer is.
\end{Dfn}
\begin{Pro} Assume that $E$ is smooth. Let $r>1$. Then the $L^r$-distortion
function
$D^X_{n,r}$ is Gateaux-differentiable at every admissible $n$-tuple
$(a_1,\ldots,a_n)$ with a Gateaux differential given by
\[
\nabla D^X_{n,r}(a_1,\ldots,a_n) = r\left(\E \left(\mbox{\bf
1}_{C_i(\alpha)\setminus\{a_i\}}(X)
\|X-a_i\|^{r-1}\nabla\|\,.\,\|(a_i-X)\right)\right)_{1\le i\le n} \in (E^*)^n
\]
where $\{ C_i (\alpha) : 1 \le i \le n \}$ denotes any Voronoi partition
induced by $\alpha = \{ a_1 , \ldots , a_n \}$.
If the norm is Fr\'echet-differentiable at every $x\neq 0$, then $\nabla
D^X_{n,r}(a_1,\ldots,a_n) $ is the Fr\'echet derivative. Furthermore, if $E$
is uniformly smooth, then $(a_1,\ldots,a_n)\mapsto \nabla
D^X_{n,r}(a_1,\ldots,a_n) $ is continuous on the set of admissible $n$-tuples
(where
$E^*$ is endowed with its norm).

When $r=1$, the above results extend to admissible $n$-tuples with
$\P_{_X}(\{a_1,\ldots,a_n\})=0$.
\end{Pro}

\noindent {\bf Remark.} In case $E =L^1$, the above proposition as well as
Proposition 1(b) do not apply since the $\|\,.\,\|_{_1}$-norm is neither
smooth nor strictly convex.

\bigskip
\noindent {\bf Proof.} A
straightforward adaptation of Lemma 4.10 in \cite{GRLU} yields both
differentiability properties. Then, if
$E$ is uniformly smooth, the mapping $x\mapsto \nabla\|\,.\,\|(x)$ is
continuous (see [2]). One derives the continuity of
$\nabla D^X_{n,r}$ by the Lebesgue dominated convergence theorem using that $
\nabla \|\,.\,\|$ takes its values in the unit ball of
$E^*$. \hfill{$\Box$} 

\begin{Dfn} Let $E$ be a   Banach space and let $r\ge 1$. A $n$-quantizer $a =
\{ a_1 , \ldots , \alpha_n \} \subset E$ of size $n$ is called
$L^r$-stationary for $X$ if $\P_X (C_i(\alpha)) > 0$ and
%
\begin{equation}\label{StatioDef}
\E \left(\mbox{\bf
1}_{C_i(\alpha)\setminus\{a_i\}}(X) \|X-a_i\|^{r-1}\nabla\|\,.\,\|(a_i -
X)\right)=0,\qquad i=1,\ldots,n,
\end{equation}
where $\{ C_i (\alpha) : 1 \le i \le n \}$ denotes any Voronoi partition
induced by $\alpha$. (This requires that
the Gateaux-differential $\nabla\|\,.\,\|(a_i-x) $ is defined $\P_{_X}(d
x)$-$a.e.$ on $C_i(\alpha)\setminus\{a_i\}$ and, furthermore, that
$\P(X\in\alpha)=0$ when $r=1$).
\end{Dfn}

This finally leads to the following proposition which makes the (expected)
connection between optimality and stationarity.
\begin{Pro}\label{OptitoStatio} Assume that $E$ is smooth and strictly convex.
Let $r>1$.  Assume that
card$({\rm supp}\,\P_{_X} ) \ge n$. Then any $L^r$-optimal
$n$-quantizer $\alpha$ is $L^r$-stationary (and admissible) for $X$. This
extends to $r=1$ if $\P_{_X}(\alpha)=0$.
\end{Pro}

\noindent {\bf Proof.} Any  $L^r$-optimal
$n$-quantizer $\alpha = \{ a_1 , \ldots , a_n \}$ is admissible by
Proposition~\ref{Pro2.1&2}$(b)$, hence the Gateaux-differential
$\nabla D^X_{n,r}(a_1 , \ldots , a_n)$ does exist and is $0$ which exactly
means stationarity. \hfill{$\Box$}
\subsection{Stationarity  for stochastic processes} Let $(T, {\cal B},\mu)$ be
a finite measure space, let
$X =(X_t)_{t\in T}$ be a bi-measurable $\R^d$-valued process defined on a
probability space $(\Omega,{\cal A},\P)$ and let $p,\,r\!\in[1,+\infty)$.
Assume  that $L^p_{\R^d} (\mu)$ is separable and that
 $\|X\|_{_p}\!\in L^r(\P)$ i.e.
\begin{equation}\label{DefXLp}
\E \left(\int_T |X_t|^p_p d \mu (t) \right)^{r/p}<+\infty.
\end{equation}

Then, the process $X$ can be seen as a (Radon) random vector taking  its
values in the Banach space $(E,\|\,.\,\|)=(L^p_{\R^d}(\mu) ,\|\,.\,\|_{_p})$
satisfying an $L^r$-integrability property, that is
$X\!\in L^r_{L^p_{\R^d}}  (\P)$. When $p\neq 1$, the $L^p_{\R^d}$-spaces are
uniformly smooth and strictly convex, so the above abstract results apply.
Furthermore, if $q$ denotes the conjugate
H\"older exponent of $p$, for every $f = (f_1 , \ldots , f_d) \in L^p_{\R^d}$,
$f\not \equiv 0$,
\[
\nabla\|\,.\,\|_{_p}(f)= \left(
\left(\frac{|f_j|}{\|f\|_{_p}}\right)^{p-1}\!\!\!\!\!{\rm sign} f_j \right)_{1
\leq j \leq d} \! \in E^{*} = L^{q}_{\R^d}
\]
so that the
$(L^r,\|\,.\,\|_{_p})$-stationarity condition reads for any Voronoi partition
$\{ C_i (\alpha) : 1 \leq i \leq n \}$ with $\P_X (C_i(\alpha)) > 0$, for
every $i$,
\begin{equation}\label{StationnaireLp}
\E\left(\mbox{\bf 1}_{C_i(\alpha)}(X)\|X-a_i\|_{_p}^{r-p}|a_{ij} -X_j
|^{p-1}{\rm sign}(a_{ij}-X_j )\right) \stackrel{L^q}{=}  0,\qquad
i=1,\ldots,n, j = j , \ldots , d
\end{equation}
with the convention $\frac{0}{\| 0 \| }=0$, where
$a_i = ( a_{i1} , \ldots a_{id})$. When $p=1$, the condition is formally the
same. This may be written in a more synthetic way by introducing the
$\alpha$-quantization $\widehat X:= \widehat X^\alpha$ of
$X$ defined by~(1.3), namely:
\begin{equation}\label{StationnaireLpbis}
\E\left(\|X-\widehat X\|_{_p}^{r-p}|X_j -\widehat X_j |^{p-1}{\rm
sign}(\widehat X_j -X_j)\, | \,\widehat X\right) \stackrel{L^q}{=} 0.
\end{equation}
When  $p=2$, $r\ge 2$ (and $\P(X \in \alpha)=0$ if $r>2$),
Equation~(3.3) looks
simpler and reads
\begin{equation}\label{quadrastatio}
a_i  \stackrel{L^2_{\R^d}}{=} \frac{ \E(X\mbox{\bf
1}_{C_i(\alpha)}(X)\|X-a_i\|^{r-2}_{_2})}{ \E(\mbox{\bf
1}_{C_i(\alpha)}(X)\|X-a_i\|^{r-2}_{_2})} , \qquad 1\le i\le n.
\end{equation}

One derives from Proposition 7 and Proposition 1 the following corollary.
\begin{Cor} Let $p,r\!\in[1,+\infty)$, let $n\ge 1$. If\begin{equation}
\left\{\begin{array}{lcl}
p,\, r>1 &\mbox{and}& \mbox{card}({\rm supp} \P_{_X})\ge n,\\
*[.4em] p>1,\; r=1
  &\mbox{and}&  \P_{_X} \mbox{ is continuous,}\\
*[.4em]  p=1,\, r\ge 1 &\mbox{and}& \P_{_{X_{j,t}}} \;\mbox{ is }\; \mu (dt)
\mbox{-}a.e.\mbox{ continuous for every} \; j \in \{ 1, \ldots , d \} ,
\end{array}\right.
\end{equation}
then, any $(L^r,\|\,.\,\|_{_p})$-optimal $n$-quantizer is
$(L^r,\|\,.\,\|_{_p})$-stationary in the sense
of (3.3).
\end{Cor}
{\bf Proof.} It remains to consider the case $p=1$. The space $L^1_{\R^d}$ is
not smooth. However, $\| . \|_1$ is Gateaux-differentiable at every $f$ such
that $f_j(t) \not= 0$ $\mu (dt)$ - a.e. for every $j$.
Now, by the Fubini Theorem, one has for every $g \in L^1(\mu)$
\[
\int_\Omega  \mu (t : X_{j,t} (\omega) = g (t)) \P (d \omega) = \int_T
\P(X_{j,t} = g (t)) \mu (dt) = 0
\]
i.e. $(X_{j,t} - g(t) \not= 0\; \mu(dt)$-a.e.) $\P$-a.s. Let $\alpha = \{ a_1
, \ldots , a_n\}$ be an
$(L^r , \| \cdot \|_1)$-optimal $n$-quantizer
and $\P_i := \P ( \cdot | \{ X \in C_i(\alpha)\})$. This definition is
consistent since
$\P(X \in C_i(\alpha)) > 0$ by Proposition 1 (a). It follows easily that
$\Psi_i : f \mapsto \int \| X-f \|^r_1 d \P_i$, $f \in L1_{\R^d}$,
is Gateaux differentiable with a Gateaux-differential given by

\[
\bigtriangledown \Psi_i(f) = \left( r \int \| X-f \|^{r-1}_1 \; \mbox{sign}
\;( f_j - X_j)d \P_i \right)_{1 \leq j \leq d} \in L^\infty_{\R^d} .
\]
Now, still following Proposition $1(a), a_i$ is a minimum for $\Psi_i$ so that
its Gateaux differential is zero. Hence, for every
$i \in \{ 1 , \ldots , n \}, j \in \{ 1, \ldots , d \}$,
\[
\int 1_{C_i (\alpha)} (X) \| X- a_i \|^{r-1}_1 \; \mbox{sign} (a_{ij} - X_j) d
\P = 0 .
\]
\hspace*{\fill}{$\Box$}

\bigskip
\noindent  {\bf Remark.} Continuity of $\P_{X_j,t} \; \mu (dt)$ - a.e. for
some $j$ implies continuity of $\P_X$.
\subsection{Pathwise regularity of stationary quantizers ($1\le p\le
r<+\infty$)} 
As before, let $E = L^p_{\R^d}(\mu)$ for some finite measure space $(T, {\cal B}, \mu)$ such that $E$ is separable.
We will derive from Equations (3.3) (and ~(3.5))  some pathwise
continuity result for the
$(L^r,\|\,.\,\|_{_p})$-stationary quantizers (which extends a result
established in \cite{LUPA1} in the purely quadratic case $p=r=2$).
For $q \in (0, \infty)$, if $X_t \in L^q_{\R^d}( \P)$ for every $t \in T$,
define the {\em ``intrinsic" semimetric}
$\rho^q_X$ on $T$ by
\[
\rho^q_X (s,t)  :=  ( \E | X_s - X_t |^q_q)^{1/(q \vee 1)}
  =  \| X_s - X_t \|^{q / (q \vee 1) }_{L^q_{\R^d} (\P)}  , \;s, t \in T.
\]
{{\large \bf Theorem 5}}
{\em Let $p,r\!\in [1,+\infty)$, $r\ge p$. Let $X$ be a bi-measurable
$\R^d$-valued process satisfying (3.2) and
\[
\forall\, t\!\in T,\qquad X_t\!\in L^{r-1}_{\R^d} (\P).
\]
Let $\alpha = \{ a_1,\ldots,a_n\}$  be an
$(L^r,\|\,.\,\|_{_p})$-stationary
$n$-quantizer (in the sense of~(3.3)). 
Set $I_r (\alpha) := \{ i \in \{ 1 , \ldots , n\} : \; \P (X = a_i) = 0 \} $ if
$r > p$ and $I_r (\alpha) := \{ 1, \ldots , n \}$ otherwise. \\
\smallskip
\noindent (a) Let $T$ be a compact metric space and let $\mu$ be a continuous
finite Borel measure on $T$. If $p=1$, if $X$ is pathwise continuous
with
\[
\mbox{supp} ( \P_X) = \{ f \in C_{\R^d} (T) : f (t) = x, t \in T_0 \}
\]
(in case $X$ is viewed as a $(C_{\R^d} (T), \| \cdot \|_{\mbox{sup}} )$-random
vector) for some $x \in \R^d$ and some closed subset $T_0 $ of $T$ with
$\mu(T_0) = 0$ and if the distribution $\P_{X_{j,t}}$ is continuous on $\R$
for every $t \in T \setminus T_0, j \in \{ 1 , \ldots , d\}$, then the
components of
$\alpha$ have $\mu$-versions consisting of continuous functions such that
$a_i (t) = x , t \in T_0 , i = 1 , \ldots , n $. \\
\smallskip
\noindent (b) If $p \in (1, \infty)$, then the components $a_i , i \in I_r ( \alpha)$
of $\alpha$ have $\mu$-versions consisting of
$\rho^{r-1}_X$-continuous functions. Furthermore, if $X_t = x \in \R^d, t \in
T_0 \subset T$, then there are such versions with
$ a_i (t) = x , t \in T_0$. \\
\smallskip
\noindent (c) If $p = 2$, then the components $ a_i , i \in I_r (\alpha)$ of $\alpha$ have $\mu$-versions
consisting of $\rho^{r-1}_X$-Lipschitz continuous functions.} \\ \\
\noindent
{\bf Remarks.} $\bullet$ If $\P_X$ is continuous then $I_r (\alpha) = \{ 1 , \ldots ,n \}$. \\
$\bullet$ If $r \geq p = 2, \E X = 0$ and $\E \| X \|^{2r-4}_2
< \infty$, then $\{ a_i : i \in I_r (\alpha) \}$ even lies in the reproducing
kernel Hilbert space of $X$. This is
a consequence of (3.5).\\
$\bullet$ Let $(T, \rho)$ be a separable metric space and $\mu$ a finite Borel measure
on $(T, \rho)$. If $p > 1$ and $t \mapsto X_t$ from $(T, \rho)$ into $L^{r-1}_{\R^d} (\P)$ is continuous that is
$\rho^{r-1}_X$ is majorized by the initial metric $\rho$ on $T$,
then the $I_r(\alpha)$-components of $\alpha$ have versions consisting of $\rho$-continuous
functions. The $L^{r-1}_{\R^d} (\P)$-continuity assumption is
fulfilled e.g. if $X$ is pathwise $\rho$-continuous
and $\| X \|_{\sup} \in L^{r-1} (\P)$. \\ \\
{\bf Proof of Theorem 5.} For every $i\!\in I_r (\alpha)$, set $
\Q_{i,r} = \mbox{\bf 1}_{C_i(\alpha)}(X)\|X - a_i\|_{_p}^{r-p}.\P$.
The measure $\Q_{i,r}$ is finite: if $r=p$, this is obvious, otherwise,
\[
\Q_{i,r}(\Omega) \le  \E\,\|X - a_i \|_{_p}^{r-p} \le \left(\E\|X -
a_i\|^r_{_p}\right)^{1-\frac pr}<+\infty.
\]
On the other hand, $\Q_{i,r}$ is a nonzero measure equivalent to
$1_{C_i(\alpha)} (X). \; \P$ since
$\P(X\in C_i(\alpha))>0$ and for $r > p$, $\P(X =a_i)=0$. Now, define on $\R
\times T$ the function $\Phi_{i j} $ by
\[
\Phi_{i j}(y,t):= \int_\Omega \varphi_{p-1}(y-X_{j,t}) d\Q_{i,r}\qquad \mbox{
where }\qquad \varphi_q(x) = \mbox{\rm sign}(x)|x|^{q}.
\]
First note that the function $\Phi_{ij}$
is real valued. If $r>p>1$, the Young inequality with $p'= \frac{r-1}{p-1}$
and $q'= \frac{r-1}{r-p}$ implies
\begin{eqnarray*}
|y-X_{j,t} |^{p-1}\|a_i-X\|_{_p}^{r-p}&\le& C(|y-X_{j,t}|^{r-1}+
\|X-a_i\|_{_p}^{r-1})\\
&\le & C (|y|^{r-1}+\|a_i\|_{_p}^{r-1}+|X_{j,t}|^{r-1}+\|X\|_{_p}^{r-1})
\end{eqnarray*}
so that $|y-X_{j,t} |^{p-1}\|a_i-X\|_{_p}^{r-p} \!\in L1(\P)$. When $r=p$ (or
$p=1$), the
result is obvious.

\smallskip
\noindent $(b)$  For every fixed $t\!\in T$ and $p>1$, $y \mapsto
\varphi_{p-1}(y-X_{j,t} )$ is (strictly) increasing, hence
$y \mapsto \Phi_{ij}(y,t)$ is strictly increasing too.  The continuity of $y
\mapsto \Phi_{ij}(y,t)$ on $\R$ for every $t \in T$ follows from the Lebesgue
dominated
convergence Theorem.  Furthermore, for every $t \in T, y \ge 0$,
\[
\Phi_{ij}(y,t)\ge \int_{\{X_{j, t} \le y \}} \varphi_{p-1}(y -
X_{j,t})d\Q_{i,r}-\int |X_{j,t}|^{p-1}d\Q_{i,r}
\]
so that $\displaystyle \lim_{y \to +\infty} \Phi_{ij} (y,t)=+\infty$ by
Fatou's Lemma. Similarly,
$\displaystyle\lim_{y \to -\infty} \Phi_{ij}(y,t)=-\infty$.

The proof reduces to providing an argument for the $\rho^{r-1}_X$-continuity
of $t\mapsto \Phi_{ij} (y,t)$ for every $y \in \R$.

\smallskip
   If $1<p\le 2$, one starts from the inequality
\[
|\varphi_{p-1}(u)-\varphi_{p-1}(v)|\le 2^{2-p}|u-v|^{p-1}\quad u,\, v\!\in\R.
\]
When  $r>p$, the H{\"o}lder inequality applied with the conjugate exponents
$\frac{r-1}{p-1}$ and $\frac{r-1}{r-p}$ yields
\[
\begin{array}{lcl}
|\Phi_{ij}(y,t)-\Phi_{ij}(y,s)|& \le & 2^{2-p}\|X_{j,t}-X_{j,s}
\|_{L^{r-1}(\P)}^{p-1}\|\,\|X- a_i\|_{_p}\|^{r-p}_{L^{r-1}(\P)} \\
& \le & 2^{2-p} ( \rho^{r-1}_X (s,t))^{\frac{p-1}{(r \wedge 2) -1}} \| \; \| X
- a_i \|_p \|^{r-p}_{L^{r-1}(\P)} .
\end{array}
\]
This still holds if $r=p$.

\smallskip
 If $p> 2$, one starts from
\[
|\varphi_{p-1}(u)-\varphi_{p-1}(v)|\le (p-1) (|u|\vee |v|)^{p-2}|u-v|,\quad
u,\, v\!\in\R.
\]
Since  $r>2$   the Holder Inequality applied with $r-1$ and $\frac{r-1}{r-2}$
yields
\begin{eqnarray*}
|\Phi_{ij}(y,t)\!-\!\Phi_{ij}(y,s)|\!\!&\!\!\le \!\!&\!\! (p\!-\!1)
\E\left(|X_{j,t} -X_{j,s} |\left(|y-X_{j,t}|\!\vee\! |y-X_{j,s} |\right)^{p-2}
\|X-a_i\|_{_p}^{r-p}\mbox{\bf 1}_{C_i(\alpha)}(X)\right)\\
\!\!&\!\!\le \!\!&\!\! (p\!-\!1) \|X_{j,t} \!-\!X_{j,s}
\|_{L^{r-1}(\P)}\!\!\left[\!  \E \left(\! \left(|y \!-\!X_{j,t} |\!\vee\!
|y \!-\!X_{j,s}
|\right)^{\frac{(p-2)(r-1)}{r-2}}\!\|X\!-\!a_i\|^{\frac{(r-p)(r-1)}{r-2}}_{_p}
\right)\! \right]^{\frac{r-2}{r-1}}.
\end{eqnarray*}
A new application of the Holder Inequality to the expectation in the right
hand side of the
above inequality yields
\begin{eqnarray*}
|\Phi_{ij}(y,t)-\Phi_{ij}(y,s)|\!\!&\!\!\le \!\!&\!\!
(p-1)\|X_{j,t}-X_{j,s}\|_{L^{r-1}(\P)}\||y-X_{j,t}|\vee
|y-X_{j,s} |\|_{L^{r-1}(\P)}^{p-2}
\left\|\,\|X-a_i\|_{_p}\right\|_{L^{r-1}(\P)}^{r-p}\\
*[.5em]\!\!&\!\!\le \!\!&\!\! C_{p,a_i} \|X_{j,t} - X_{j,s}
\|_{L^{r-1}(\P)}\left(|y|^{p-2}+\|X_{j,s} \|^{p-2}_{L^{r-1}(\P)}+\|X_{j,t}
\|^{p-2}_{L^{r-1}(\P)}\right)\\
& \le & C_{p, a_i} \rho^{r-1}_X (s,t) (| y |^{p-2} + \| X_s
\|^{p-2}_{L^{r-1}_{\R^d}(\P)} + \| X_t \|^{p-2}_{L^{r-1}_{\R^d}(\P)} ).
\end{eqnarray*}

Owing to these properties, one easily checks that for every $t\!\in T$, the
equation $\Phi_{ij} (y,t)=0$ admits a unique solution
$y_{ij} (t)$  and that
the implicitly defined function $t\mapsto y_{ij}(t)$ is $
\rho^{r-1}_X$-continuous. On the other hand the function $a_i$
satisfies $\mu(dt)$-a.e. $\Phi_{ij} (a_{ij} (t),t)=0$ so that $y_{ij} (t) =
a_{ij}(t)$
$\mu(dt)$-a.e..

If $X_t = x \in \R^d, t \in T_0$ then $\Phi_{ij}(y,t) = \varphi_{p-1} (y -
x_j) \Q_{i,r} ( \Omega), t \in T_0$ so that $y_{ij} (t) = x_j$.

\smallskip
\noindent $(a)$  Now let $T$ be a compact metric space. When $p=1$,
\[
\Phi_{ij}(y , t) = \int_\Omega {\rm sign }(y-X_{j,t} )\, d\Q_{i,r}.
\]
The continuity of $y \mapsto \Phi_{ij} (y,t)$ on $\R$ for every $t \in T
\setminus T_0$ and the continuity of $t \mapsto \Phi_{ij}(y,t)$
at every point $t \in T \setminus T_0$ for every $y \in \R$
follows from   the pathwise continuity of $X$ and from the continuity  of
$\P_{X_{j,t}}, t \in T \setminus T_0$, by  the Lebesgue dominated convergence
Theorem: the sign function is bounded and $\Q_{i,r} \ll \P$. Similarly
one shows that
$\displaystyle
\lim_{y \to \pm \infty} \Phi_{ij}(y,t)=\pm \Q_{i,r}(\Omega) \; \forall t \in T
$.
It is also obvious that $y \mapsto \Phi_{ij}(y , t)$ is nondecreasing $\forall
t \in T$. To establish strict monotonicity
$\forall t \in T \setminus T_0$, one proceeds as follows: let us consider the
subset of $C_i(\alpha)$ defined by
\[
U_i(\alpha) := \{ f \in C_{\R^d} (T,T_0) : \| f - a_i \|_1 < \min_{j \not= i}
\| f - a_j \|_1 \}
\]
where
\[
C_{\R^d} (T,T_0) := \{ f \in C_{\R^d} (T) : f(t) = x, t \in T_0 \} .
\]
It is a nonempty open subset of $(C_{\R^d}(T,T_0) , \| \cdot \|_1)$ since
$C_{\R^d}(T,T_0)$ is everywhere
$\| \cdot \|_1$-dense in $L^1_{\R^d}(\mu)$ in view of $\mu(T_0) = 0$. Now, for
$t \in T \setminus T_0$ and every nonempty open interval $I$
the set $\{
f \in C_{\R^d}(T,T_0) : f_j(t) \in I \}$ is clearly everywhere dense in
$(C_{\R^d}(T,T_0), \|  \cdot \|_1)$ since $\mu(\{t\}) = 0$ so that
\[
U_i(\alpha) \cap \{ f \in C_{\R^d}(T,T_0) : f_j(t) \in I \}
\]
is a nonempty set. On the other hand, $f \mapsto \| f \|_1$ and $f \mapsto
f_j(t)$ are both continuous as functionals on
$(C_{\R^d}(T,T_0), \| \cdot \|_{_{\sup}})$ so that
$U_i(\alpha) \cap \{ f \in C_{\R^d} (T,T_0):f_j (t) \in I \}$ is a (nonempty)
open subset of
$(C_{\R^d} (T,T_0), \|  \cdot \|_{\sup} )$.
Now, if
$\Phi_{ij} (y,t)= \Phi_{ij} (y',t)$ for some $y < y'$, then $\Q_{i,r}(X_{j,t}
\!\in(y , y'))=0$.
$\Q_{i,r}$ is equivalent to
$\mbox{\bf 1}_{C_i(\alpha)}(X).\P$. Consequently
\[
\P (\{X\!\in U_i(\alpha)\}\cap  \{X_{j,t} \!\in (y ,y')\})=0.
\]
This is impossible owing to the assumption on the support of $\P_{_X}$.
Consequently $y\mapsto \Phi_{ij}(y,t)$ is strictly increasing for every $t \in
T \setminus T_0$ and one
concludes like in the case $p>1$ to the existence of a continuous version of
$\alpha$ in $C_{\R^d}(T,T_0)$.

To be a bit more precise, the equation $\Phi_{ij} (y,t) = 0$ has for $t \in T
\setminus T_0$ a unique solution $y_{ij}(t) \in \R$ and for $t \in T_0$, since
$X_{j,t} = x_j \; \P$-a.s., $y_{ij} (t) = x_j$ is the unique solution. The
function $y_{ij} : T \rightarrow \R$ is continuous at every $t \in T \setminus
T_0$ since $\Phi_{ij}( \cdot , t)$ is strictly increasing on $\R$ and
$\Phi_{ij} (y, \cdot)$ is continuous at $t$ for every $y \in \R$. One must
consider the behaviour of $y_{ij}$ at $t \in T_0$ more carefully. First note
that $\Phi_{ij} (y, \cdot)$ is continuous at $t \in T_0$ for every $y \not=
x_j$ since $X_j$ is pathwise continous.
Now let $(s_n)$ be a sequence in $T$ going to $t$ such that $y_{ij}(s_n) \geq
x_j + \eta$ for some $\eta > 0$.
Then, $\Phi_{ij}(x_j + \eta, s_n) \leq \Phi_{ij} (y_{ij} (s_n) , s_n) = 0$ for
every $n \geq 1$ so that
$\mbox{sign}(\eta) \Q_{i,r} (\Omega) = \Phi_{ij} (x_j + \eta , t) = \lim_{n
\to \infty} \Phi_{ij} (x_j + \eta , s_n) \leq 0$ which is
impossible. Hence $\limsup_{s \to t} y_{ij} (s) \leq x_j$. One shows similarly
that $\liminf_{s \to t} y_{ij} (s) \geq x_j$ i.e.
$\lim_{s \to t} y_{ij} (s) = x_j = y_{ij} (t)$. \\
\smallskip
\noindent $(c)$ It is a   consequence of Equation (3.5):
\[
|a_{ij}(t)-a_{ij}(s)| = \frac{\E( | X_{j,t} -X_{j,s} | L_i)}{\E( L_i)}\quad
\mbox{ with } \quad L_i= \mbox{\bf 1}_{C_i(\alpha)} (X)
\|X-a_i\|^{r-2}_{_2}.
\]
When $r>2$, The Holder Inequality yields the announced result
\[
\max_{i \in I_r (\alpha)} | a_i (t) - a_i(s) |_{_{1}} \leq C_{X,\alpha}
\rho^{r-1}_X (s,t)
\]
with
\[
C_{X,\alpha}:= d \max_{i \in I_r (\alpha) } (\E( \|X-
a_i\|^{r-1}_{_2}))^{(r-2)/(r-1)}/(\E(\mbox{\bf 1}_{C_i(\alpha)} (X) \|X -
a_i\|^{r-2}_{_2})).
\]
When $r=2$, one sets accordingly $C_{X,\alpha}:=1/\min_{1\le i\le n} \P(X\!\in
C_i(\alpha)) $. \hfill{$\Box$} \\ \\
%
{\sc Examples} First consider real or $\R^d$-valued processes with
$T = [0,t_0]$ and $\mu(dt) = dt$. \\
\smallskip
\noindent $\bullet$ The $(L^r,\|\,.\,\|_{_p})$-stationary $n$-quantizers,
$1\le p\le r<+\infty$, of the {\em standard Brownian motion} and are made up
with continuous functions which are null at $0$, $1/2$-H{\"o}lder if $p=2$.
The same result holds for the {\em Brownian bridge} over $[0,t_0]$ where
any of its  stationary quantizers are $0$ at $t_0$ and for the {\em standard
$d$-dimensional Brownian motion}. \\
\smallskip
\noindent $\bullet$ One considers a $\R^d$-valued {\em Brownian diffusion
process}
\[
\begin{array}{lcl}
d X_t & = & b(t, X_t) dt + \sigma (t, X_t) d W_t, t \in [0, t_0] \\
X_0 & = & x, x \in \R^d ,
\end{array}
\]
where $W$ is a $m$-dimensional standard Brownian motion and $b : [0,t_0]
\times \R^d\rightarrow \R^d$,
$\sigma : [0, t_0] \times \R^d \rightarrow \R^{d \times m}$ are Borel
functions with linear growth such that the above SDE admits at least one
(weak) solution over $[0,t_0]$.
This
 solution is pathwise continuous and it is classical background
(see~\cite{KASH}) that
$ \| X \|_{\sup} \in L^r(\P)$ for every $r \in (0, \infty)$ and
\[
\E | X_s - X_t |^q_q \leq C_q | s-t | ^{q/2}
\]
for every $q \in (0, \infty)$. Thus the $(L^r, \| \cdot \|_p)$-stationary
$n$-quantizers, $1 < p = r < \infty$, are made up with continuous functions
which are
$x$ at $t = 0$, $1/2$-Holder if $p = 2$. The same holds if $1 \leq p \leq r <
\infty$ for the
homogeneous SDE with $b$ and $\sigma$ independent of $t$ and $d = m$ provided
$b_i $ and $\sigma_{ij}$ are bounded with bounded derivatives up to order 3
and $\sigma \sigma^T$ is uniformly elliptic. In fact, the assumptions imply
that $\P_{X_{j,t}}$ has a Lebesgue density for every
$j \in \{ 1, \ldots , d\}, t \in (0,t_0]$
and by the support theorem, in $C_{\R^d} ([0,t_0])$,
\[
\mbox{supp} (\P_X) = \{ f \in C_{\R^d} ([0,t_0]) : f (0) = x \}
\]
(see [3], p. 11 and [1], p. 25). \\
$\bullet$ The {\em fractional Brownian motion} $W^H$ on $[0, t_0]$ with Hurst
exponent $H \in (0,1)$ is a centered continuous Gaussian process having the
covariance function
\[
E W^H_s W^H_t = \frac{1}{2} ( | s |^{2H} + | t |^{2H} - | s-t |^{2H} )
\]
and thus satisfies for every $q \in (0, \infty)$
\[
\E | W^H_s - W^H_t |^q = C_{H,q} | s-t |^{qH} .
\]
Consequently, $(L^r, \| \cdot \|_p)$-stationary $n$-quantizers, $1 \leq p \leq
r < \infty$
are made up with continuous functions which are null at $t = 0$, H-H{\"o}lder
if $p=2$.

\smallskip
\noindent $\bullet$ We consider some examples of (cadlag) real L\'evy
processes $X = (X_t)_{t \in \R_{+}}$ restricted to
$[0, t_0]$ (without Brownian component). Since the increments of $X$ are
stationary and $X_0 = 0$,
\[
\E | X_s - X_t |^q = \E | X_{| s-t| } |^q
\]
so that the behaviour of the semimetric $\rho^q_X$ reduces to the behaviour of
$t \mapsto \E | X_t |^q$.

- The {\em $\rho$-stable L\'evy motions} indexed by $\rho \in (0,2)$ satisfy a
self-similarity property, namely
\[
X_t \stackrel{d}{=} t^{1/\rho} X_1 .
\]
Furthermore,
\[
\sup \{ q > 0: \E | X_1 |^q < \infty \} = \rho \; \mbox{and} \; \E | X_1
|^\rho = \infty .
\]
For this background see \cite{SAMORO}. It follows that for every $q \in (0,
\rho)$
\[
\E | X_t |^q = t^{q/\rho} \E | X_1 |^q < \infty .
\]
Consequently, since the $\rho$-stable distributions $\P_{X_t}, t > 0$ have a
Lebesgue density, the
$(L^r, \| \cdot \|_p)$-stationary $n$-quantizers, $1 < p \leq r < \rho$,
are made up with continuous functions which are null at 0.

- The {\em $\Gamma$-processes} are L\'evy processes whose distribution
$\P_{X_{t}}$ at $t > 0$ is a $\Gamma (a,t)$-distribution
\[
\P_{X_t}(dx) = \frac{a^t}{\Gamma(t)} 1_{(0, \infty )} (x) x^{t-1} e^{-ax} dx ,
\]
$a > 0$. So, for every $q > 0$
\[
\E | X_t |^q = \frac{\Gamma(t+q)}{a^q\Gamma(t+1)} t .
\]
Consequently, $(L^r, \|  \cdot \|_p)$-stationary $n$-quantizers, $1 < p \leq
r < \infty$, are made up with continuous functions, $1/(r-1)$-H\"older if $p =
2$.

- The {\em compound Poisson process} is given by $X_t = \sum\limits^{N_t}_{j=1} U_j$, where
$U_1, U_2 , \ldots$ are i.i.d. real random variables with $\P(U_1 = 0) = 0$
and $N =(N_t)_{t \geq 0}$ is a standard Poisson process
(with intensity $\lambda$) independent of $(U_j)_{j \geq 1}$. If $q \in (0, \infty)$ and 
$\E \mid U_1 \mid^q < \infty$, easy computations show that
\[
\E \mid X_t \mid^q \leq \E \mid U_1 \mid^q \E N^{1 \vee q}_t \leq C_{q, \lambda, U} t < \infty
\]
Assume $\E \mid U_1 \mid^r < \infty$. Then the $(L^r, \parallel \cdot \parallel_p)$-stationary $n$-quantizers, $1 < p \leq r < \infty$, are made up with continuous functions, $1/(r-1)$-H\"older if $p = 2$. 
Here it has to be noticed that the function $f = 0$ is the only atom of
$\P_X$ in $L^p([0, t_0], dt)$.

Theorem 5$(a)$ does not apply to the above examples because of
the pathwise continuity assumption so that the  case $p=1$ remains open.

As for a real multiparameter process on $T = [0,t_0]^k$ with $\mu (dt) = dt$:
\\
$\bullet$ The $(L^r, \| \cdot \|_p)$-stationary $n$-quantizers, $1 \leq p \leq
r < \infty$, of the {\em standard Brownian sheet} are made up with continous
functions which are null on $\bigcup\limits^{k}_{i=1} \{ t \in T: t_i = 0 \}$
and 1/2-H{\"o}lder if $p=2$. \\ 

As for an example with noncompact $T$ consider $T = \R_{+}$ and $\mu (dt) = e^{-bt}dt, b > 0$. \\
$\bullet$ The stationary {\em Ornstein-Uhlenbeck process} $X = (X_t)_{t \geq 0}$ on $\R_{+}$ is a centered continuous Gaussian process having the covariance function
\[
\E X_s X_t = e^{-c\mid s-t \mid} , c > 0 .
\]
Clearly, $X$ can be seen as $L^p (\R_{+}, \mu)$-valued random vector for every $p \in [1, \infty)$. The process satisfies for every $q \in (0, \infty)$
\[
\E \mid X_s - X_t \mid^q = C_q (1-e^{-c\mid s-t\mid })^{q/2} .
\]
Consequently, $(L^r, \parallel \cdot \parallel_p)$-stationary $n$-quantizers, $1 < p \leq r < \infty$ have components consisting of continuous functions, $1/2 $-H\"older if $p = 2$. \\ \\
{\sc A counterexample when $p=+\infty$} We will exhibit a bounded
pathwise continuous process $X$
on $T = [0,1]$ having a discontinuous $(L^r,\|\,.\,\|_{_\infty})$-optimal
$1$-quantizer.
Consider functions $f_n \in C([0,1]), n \in \N$ and $\P_X$ from the
$C([0,1])$-counterexample following Theorem 4. Then set
$h := \frac{1}{2}(\mbox{\bf 1}_{[0,1/2]}-\mbox{\bf 1}_{(1/2,1]})$.
One checks that in $L^\infty ([0,1],dt)$, for every
$n\ge 1$,
\[
\| f_n - h \|_{_\infty}=1/2
\]
so that
$$
\forall\, r\!\in[1,+\infty], \qquad \left\|\,\|X- h
\|_{_\infty}\right\|_{L^r(\P)}=1/2.
$$
On the other hand,
\[
e_{1,1} (X, L^\infty) = e_{1,1} (X,C([0,1]) = 1/2 = \|  \; \| X - h
\|_{_\infty}
\|_{L^1 (\P)}
\]
by the $C([0,1])$-counterexample and Theorem 4.
Consequently, the $\|\,.\,\|_{L^r(\P)}$-norm being nondecreasing as a function
of $r$,
\[
\forall\, r\!\in[1,+\infty], \qquad e_{1,r}(X, L^\infty )= \| \;\|  X - h
\|_{_\infty} \|_{L^r(\P)} = 1/2
\]
with obvious definition of $e_{1, \infty}$. The function $h$ is an
$(L^r,\|\,.\,\|_{_\infty})$-optimal
$1$-quantizer without continuous $dt$-version of the pathwise continuous
process $X$, $1\le r\le +\infty$.

Note that $t\mapsto X_t$ from $[0,1]$ into $L^p(\P)$ is   continuous for any
$p\!\in[1,+\infty)$ since $X$ is pathwise continuous and uniformly bounded by
$1$. Consequently it follows from Theorem 5 that, as soon as
$1 < p \leq r < + \infty$, any $(L^r,\|\,.\,\|_{_p})$ optimal $n$-quantizer of
$X$ (has a $dt$-version which) consists of continuous functions. 
However,
$t\mapsto X_t$ from $[0,1]$ into $L^\infty(\P)$ is   not continuous (at
$t=1/2$), so  the pathwise regularity of an optimal
$(L^r,\|\,.\,\|_{_\infty})$-optimal $n$-quantizer of an
$L^\infty(\P)$-continuous process remains open. But
the $\|\,.\,\|_{_\infty}$-norm being nowhere Gateaux-differentiable, the very notion of 
$(L^r, \parallel \cdot \parallel_\infty)$-stationary quantizer no longer exists. So this would require
to develop a new approach. \\ \\
\end{document}